%% file: chow.tex
\documentclass[rotate]{article}
\usepackage{lscape,graphics}

\usepackage{ifthen}
\usepackage[latin1]{inputenc}
\usepackage{amsmath,amssymb,amsthm}
\usepackage[dvips]{graphicx,epsfig}
\usepackage{color} %pour xfig
\newcommand{\A}{\mathbb{A}}

\def\CC{{\mathbb C}}

\def\FF{{\mathbb F}}
\def\GG{{\mathbb G}}
\def\HH{{\mathbb H}}

\def\NN{{\mathbb N}}

\def\PP{{\mathbb P}}
\def\QQ{{\mathbb Q}}
\def\RR{{\mathbb R}}

\def\ZZ{{\mathbb Z}}
\def\fxa{{\ensuremath \mathcal A}}
\def\fxb{{\ensuremath \mathcal B}}
\def\fxc{{\ensuremath \mathcal C}}

\def\fxf{{\ensuremath \mathcal F}}

\def\fxi{{\ensuremath \mathcal I}}

\def\fxo{{\ensuremath \mathcal O}}
\def\fxp{{\ensuremath \mathcal P}}

\def\fxx{{\ensuremath \mathcal X}}

\newcommand{\fb}{\ensuremath{\downarrow}}

\newcommand{\fd}{\ensuremath{\rightarrow}}

\newcommand{\findem}{\nolinebreak\vspace{\baselineskip} \hfill\rule{2mm}{2mm}\\}
\renewcommand{\phi}{\ensuremath{\varphi}}
\newcommand{\inc}{\ensuremath{\subset}}

\newcommand{\ox}{\otimes }
\newcommand{\pla}{\A^2}
\newcommand{\plp}{\PP^2}
\renewcommand{\phi}{\ensuremath{\varphi}}
\newcommand{\x}{\ensuremath{\times}}

\newcommand{\s}{Spec\ }

\newtheorem{nt}{Notation}

\newtheorem{prop}[nt]{Proposition}

\newtheorem{exercice}{Exercice} 
\newcounter{numeroquestion} 
\newcounter{numerosousquestion}

\newcommand{\sousquestion}{\ifthenelse{\value{numerosousquestion}=1}{}{\\}\textbf{\roman{numerosousquestion})} \addtocounter{numerosousquestion}{1}}
\newtheorem{correction}{Correction}

\newenvironment{listecompacte}
{\begin{list}
    {\ensuremath{\bullet}}
    {\setlength{\topsep}{2pt}
      \setlength{\itemsep}{1pt} \setlength{\parsep}{0pt}}
}
{\end{list}
}

\newtheorem{coro}[nt]{Corollary}

\newtheorem{lm}[nt]{Lemma} 
 
\newtheorem{rem}[nt]{Remark} 
\newtheorem{thm}[nt]{Theorem}

\newenvironment{dem}{\noindent\textit{Proof.} }{\findem}

\begin{document}
\sloppy
\title{The Chow ring of punctual Hilbert schemes on toric surfaces}
\date{}
\author{Laurent Evain (laurent.evain@univ-angers.fr)}
\maketitle
\newcommand{\esc}{{\cal{E}}}
\newcommand{\mesc}{{\cal{M}\cal{E}}}
\newcommand{\hilb}{{\cal{H}}}
\newcommand{\mhilb}{{\cal{M}\cal{H}}}
\newcommand{\Part}{{\cal{P}}art}

%%%%%%%%%%%%%%%%%%%%%%%%%%%%%%%%%%%%%%%%%%%%%%%%%%%%%%%%%%%%%%%%
%%%%                    Debut du contenu                  %%%%%
%%%%%%%%%%%%%%%%%%%%%%%%%%%%%%%%%%%%%%%%%%%%%%%%%%%%%%%%%%%%%%%%

\section*{Abstract: } Let  $X$ be a smooth projective toric surface, 
and $\HH^d(X)$ the Hilbert scheme parametrising the length $d$
zero-dimensional  
subschemes of $X$. We compute the rational Chow ring $A^*(\HH^d(X))_\QQ$. More
precisely, if $T\subset X$ is the two-dimensional torus contained in $X$,
we compute the rational equivariant Chow ring $A_T^*(\HH^d(X))_\QQ$ and the usual Chow
ring is an explicit quotient of the equivariant Chow ring. The case of some
quasi-projective toric surfaces such as the affine plane are
described by our method too.

\section*{Introduction}
\label{sec:introduction}
Let $X$ be a smooth projective surface and 
$\HH^d(X)$ the Hilbert scheme 
parametrising the zero-dimensional subschemes of
length $d$ of $X$.
The problem is to compute the  rational cohomology
$H^*(\HH^d(X))$.
 
The additive structure of the cohomology is well
understood. 
First, Ellingsrud and Str{\o}mme
computed in \cite{ellingsrud-stromme87:chow_group_of_hilbert_schemes}
the Betti
numbers  $b_i(\HH^d(X))$ when $X$ is a plane or a Hirzebruch surface
$\FF_n$. 
The Betti numbers $b_i(\HH^d(X))$ 
for a smooth surface $X$ were computed by G{\"o}ttsche
\cite{gottsche90:nbBettiDuSchemaHilbertDesSurfaces}
who realised them as coefficients of an explicit power series
in two variables depending on the Betti numbers of $X$.
This nice and surprising organisation of the Betti numbers as 
coefficients of a power series was explained by Nakajima in terms of 
a Fock space structure on the cohomology of the Hilbert
schemes \cite{nakajima97:_heisenberg_et_Hilbert_schemes}.
Grojnowski announced similar results 
\cite{grojnowski:resultatsSimilairesANakajima}. 

As to the multiplicative structure of $H^*(\HH^d(X))$, 
the picture is not as clear. There are descriptions valid for general
surfaces $X$ but quite indirect and unexplicit, 
and more explicit descriptions for some
special surfaces $X$. 

The first steps towards the multiplicative structure 
were performed again by 
Ellingsrud and Str{\o}mme 
\cite{ellingsrud-stromme93:towardsTheChowRingOfPP2}
(see also Fantechi-G\"ottsche
\cite{fantechi-gottsche93:cohomologie-3-points}) .
They gave an indirect
description of the ring structure in the case $X=\plp$ in terms
of the action of the Chern classes of the tautological bundles. 
Explicitly, Ellingsrud and Str{\o}mme constructed a variety $Y$ whose
cohomology is computable, an embedding $i:\HH^d(\plp)\fd Y$, and proved
the isomorphism $H^*(\HH^d(\plp))\simeq H^*(Y)/Ann(i_*(1))$ where
$Ann$ is the annihilator.

When $X=\pla$, Lehn gave in
\cite{lehn99:_chern_classes_of_tautological_sheaves_on_Hilbert_schemes}
an identification between the cohomology ring  $H^*(\HH^d(\pla))$ and
a ring of differential operators on a Fock space. With some
extra algebraic work, it is possible to derive from it 
a totally explicit description of the cohomology ring
$H^*(\HH^d(\pla))$. This was done by Lehn and Sorger in 
\cite{lehn_sorger01:cup_product_on_Hilbert_schemes}.
The same result was obtained independently by Vasserot 
\cite{vasserot01:anneauCohomologieHilbert}
by methods relying on equivariant cohomology. 

When $X$ is a smooth projective surface, Costello and Grojnowski have
identified $H^*(\HH^d(X))$ with 
two algebras of operators 
\cite{costello-grojnowski03:CohomoSchemaHilbertPonctuel}.

Lehn and Sorger 
extended their results to the case of $K3$ surfaces
\cite{lehn_sorger02:cup_product_on_Hilbert_schemes_for_K3}.
Li,Qin,Wang have computed the ring structure of $H^*(\HH^d(X))$ 
when $X$ is the total space of a line bundle over $\PP^1$
\cite{liQinWang04mathAG:cohomoDesSchemaHilbertSurface=FibreSurP1}.
\medskip

The goal of this work is to compute the Chow ring
$A^*(\HH^d(X))$ when $X$ is a smooth projective toric surface. 
The description is new even in the case $X=\plp$. 

Nakajima's construction 
\cite{nakajima97:_heisenberg_et_Hilbert_schemes}
has been fundamental and many of the above papers 
(\cite{vasserot01:anneauCohomologieHilbert}, 
\cite{lehn99:_chern_classes_of_tautological_sheaves_on_Hilbert_schemes},
\cite{lehn_sorger01:cup_product_on_Hilbert_schemes},
\cite{lehn_sorger02:cup_product_on_Hilbert_schemes_for_K3},
\cite{liQinWang04mathAG:cohomoDesSchemaHilbertSurface=FibreSurP1}) 
rely on it. 
The present work is independent of Nakajima's framework
and uses equivariant Chow rings as the main tool.

For simplicity, we use the notation
$\HH^d$ instead of $\HH^d(X)$. 
We use the formalism of Chow rings and work over any
algebraically closed field $k$. 
When $k=\CC$, the Chow ring co{\"\i}ncides with
usual cohomology since the action of the two-dimensional torus $T$ on 
$X$ induces an action of $T$ on $\HH^d$ with a finite number of fixed
points.  
\medskip

\textit{Equivariant Chow rings.} 
The construction of an equivariant Chow ring associated with an
algebraic space endowed with an action of a linear algebraic group has
been settled by  Edidin and Graham
\cite{edidinGraham:constructionDesChowsEquivariants}.
Their construction is modeled after the Borel in equivariant
cohomology.
Brion 
\cite{brion97:_equivariant_chow_groups}
pushed the theory further in the case the group is a torus 
$T$ acting on a variety $\fxx$. In our setting, $\fxx$ is smooth and
projective. Brion gave a description of 
Edidin and
Graham's equivariant Chow ring by generators and relations.  
This alternative construction makes it 
possible to prove that the usual Chow ring is an explicit quotient of the
equivariant Chow ring. This is the starting point
of this work: to realize the usual Chow ring as a quotient of
the equivariant Chow ring. Explicitly, the morphism $\fxx\fd \s k$
yields a pullback on the level of Chow ring and makes  $A_T^*(\fxx)$
a  $A_T^*(\s k)$-algebra. There is an
isomorphism $A^*(\fxx)\simeq A_T^*(\fxx)/A_T^{>0}(\s k)
A_T^*(\fxx)$.

Moreover, over the rationals,
the restriction to fixed points  $A_T^*(\fxx)_{\QQ}\fd A^*_T(\fxx^T)_{\QQ}$ is injective
and its image is the intersection of the images of the morphisms 
$A_T^*(\fxx^{T'})_\QQ\fd A^*_T(\fxx^T)_\QQ$ where $T'$ runs over all one
codimensional subtori of $T$. 
\\
Thus the natural context is that
of rational Chow rings and we lighten the notations: From
now on, the symbols
$A^*(\fxx), A^*_T(\fxx)$  will implicitly 
stand for the rational Chow rings 
$A^*(\fxx)_{\QQ}, A^*_T(\fxx)_{\QQ}$ . 
\\
We apply Brion's results to $\fxx=\HH^d$. The locus $\fxx ^T=\HH^{d,T}$ is a
finite number of points and the 
ring $A_T^*(\HH^{d,T})$ is a product of polynomial rings. 
In particular, the multiplicative structure of $A_T^*(\HH^{d})\subset
A_T^*(\HH^{d,T})$ is completly determined. 
In view of the above description, 
the problem of computing $A^*_T(\HH^d)\subset A^*_T(\HH^{d,T})$ 
reduces to the computation of 
$A^*_T(\HH^{d,T'})\subset A^*_T(\HH^{d,T})$. The
steps are as follows.
\begin{itemize}
\item First, we study the geometry of the locus $\HH^{d,T'}\inc
  \HH^{d}$. We identify its
  irreducible components with products $V_1\x \dots\x
  V_r$ where each term $V_i$ in the product is a projective space or a
  graded Hilbert scheme $\HH^{T',H}$, in the sense of Haiman-Sturmfels
\cite{haiman_sturmfels02:multigradedHilbertSchemes}. 
\item A graded Hilbert scheme $\HH^{T',H}$ 
  appearing as a term $V_i$ is embeddable in a product $\GG$ of
  Grassmannians. 
  A sligth modification of an argument by King and Walter
  \cite{king_walter95:generateurs_anneaux_chow_espace_modules}
  shows that the restriction morphism $A_T^*(\GG)\fd A_T^*(\HH^{T',H})$ is
  surjective. The idea for this step is that the universal family over
  $\HH^{T',H}$ is a family of $k[x,y]$-modules with a nice
  resolution. Since the equivariant Chow ring $A_T^*(\GG)$ is
  computable, we obtain a description of  $A_T^*(\HH^{T',H})$.
\item 
It then suffices to 
put the two above steps together via a K\"uneth equivariant formula
to obtain a description of the equivariant 
Chow ring $A^*_T(\HH^{d,T'})$, thus of $A^*_T(\HH^d)=\cap
_{T'}A^*_T(\HH^{d,T'})$ 
(Theorem \ref{thr:description du Chow avec
  produit tensoriel}).
\end{itemize}
At this point, the description of the equivariant Chow ring is
complete, but the formula  in theorem  \ref{thr:description du Chow avec
  produit tensoriel}
involves tensor products, direct sums and intersections. 
The last step consists in an application of a Bott formula
(proved by Edidin and Graham in an algebraic context
\cite{edidin_Graham98:formuleDeBott}) to 
get a nicer description. This is done in Theorem 
\ref{thr:description du Chow avec congruences}:
If $\hat{T}$ is the character group of $T$, $S=Sym(\hat{T}\ox \QQ)\simeq
\QQ[t_1,t_2]$ is the symmetric $\QQ$-algebra over $\hat{T}$, 
$A^*_T(\HH^d)\inc S^{\HH^{d,T}}$ is realised as a set
of tuples of polynomials satisfying explicit congruence relations. 
In this setting, the usual Chow ring is the quotient of the equivariant Chow
ring by the ideal generated by the elements $(f,\dots,f)$, $f$
homogeneous with positive degree.  We 
illustrate our method
on the case of $\HH^3(\plp)$ in 
theorem \ref{thr:leCasHilbTroisP2}.

Though we suppose for convenience that the underlying surface $X$ is projective, 
the descriptions of the Chow ring could be performed with conditions on $X$ 
weaker than projectivity: $X$ need 
only to be filtrable \cite{brion97:_equivariant_chow_groups}. 
In particular, the method applies for the affine plane.

The key results about equivariant Chow rings  used in the text
have been extended to an equivariant $K$-theory
setting 
by Vezzosi and Vistoli\cite{Vezzosi_Vistoli:KTheorieEquivarianteEtSousTores}.
Thus the method developped in the present paper should generalize to
equivariant $K$-theory as 
well.  
\medskip

\textit{Acknowledgments.} Michel Brion generously shared his
knowledge about equivariant cohomology. 
It is a pleasure to thank him for the stimulating discussions we had.

\section{Hilbert functions}
\label{sec:objectsInvolved}

To follow the method sketched in the introduction, we need to compute
the irreducible components of $\HH^{d,T'}$ for a one-codimensional
torus $T'\inc T$. 
In this section, we introduce the basic notations and several 
notions of Hilbert functions useful to describe these irreducible
components. We also define the Hilbert schemes and Grassmannians
associated with these Hilbert functions. 
 
\subsection*{The toric surface $X$}
\label{sec:toric-variety-x}

Let $T$ be a 2-dimensional torus with character group $\hat T$. 
Let $N=Hom_{\ZZ}(\hat T,\ZZ)$, $N_\RR=N\ox \RR$ and $\Delta\subset N_\RR$ a
fan defining a smooth projective toric surface $X$. Denote the maximal cones
of $\Delta$ by $\sigma_1,\dots,\sigma_r$ with the convention 
$\sigma_{r+1}=\sigma_1$, and by $p_1,\dots,p_r$ the corresponding closed
points of $X$. Assume that the cones are ordered such that
$\sigma_i\cap \sigma_{i+1}=\sigma_{i,i+1}$ is a one-dimensional
cone. Denote
respectivly by $U_{i,i+1},O_{i,i+1},V_{i,i+1}=\overline O_{i,i+1}$ 
the open subvariety, the orbit and the
closed subvariety of $X$ associated with the cone $\sigma_{i,i+1}$. 
Define
similarly $U_i\inc X$ the open subscheme associated with
$\sigma_{i}$. Explicitly, if $\sigma_i^\nu\inc \hat T\ox \RR$ is
the dual cone of $\sigma_i\inc \NN_\RR$ and 
$R_i=k[\sigma_i^\nu\cap \hat T]$,
then $U_i = \s R_i$.
The inclusion 
$R_i \inc k[\hat T]$ induces an open embedding $T=\s k[\hat
T]\hookrightarrow U_i$. The action of $T$ on itself by translation
extends to an action of $T$ on $U_i$, and to an action of $T$ on
$X=\cup U_i$. 

The open subscheme $U_i$ is isomorphic to an affine plane $\s k[x,y]$. 
When using such coordinates $x,y$, we require $xy=0$ to be the equation of 
$V_{i-1,i} \cup V_{i,i+1}$ around $p_i$. The isomorphism
$U_i\simeq \s k[x,y]$ is then defined up to the automorphism of $k[x,y]$ that
exchanges the two coordinates.

\subsection*{Subtori and their fixed locus}
\label{sec:subtori}

Let $T'\hookrightarrow T$ be a one-dimensional subtorus of $T$,
$\hat{T}'$ its 
character group. The torus $T'$ acts on $X$ by restriction 
and $U_i\inc X$ is an
invariant open subset.  The action of $T'$ on
$U_i$ induces a decomposition
$R_i=\sum_{\chi \in \hat{T}'} R_{T',i,\chi}$, where
$R_{T',i,\chi}\subset  R_i$ is the
subvector space on which $T'$ acts with character $\chi$. 
\\
One shows easily that the fixed locus
$X^{T'}$ admits two 
types of connected components. Some components are isolated fixed
points . We let  
\begin{displaymath}
PFix(T')=\{p\in X^{T'},\ p\ isolated\}.
\end{displaymath}
The other components are projective lines $V_{i,i+1}\simeq \PP^1$ joining two points
$p_i,p_{i+1}$ of $X^{T}$. 
We let
\begin{displaymath}
LFix(T')=\{ \{p_i,p_{i+1}\},\  p_i,p_{i+1} \mathrm{lie\ in\ an\
invariant\ } \PP^1\}.
\end{displaymath}
By construction,
\begin{displaymath}
  X^T=LFix(T')\cup PFix(T').
\end{displaymath}

\subsection*{Staircases and Hilbert functions}
\label{sec:staircases}
A staircase $E\subset \NN^2$ is a subset whose complement
$C=\NN^2\setminus E$ satisfies $C+\NN^2\subset C$. In our context, 
the word staircase will stand for finite staircase. 
By extension, a staircase $E\subset k[x,y]$ is a set of monomials
$m_i=x^{a_i}y^{b_i}$ such that the the set of exponents $(a_i,b_i)$
is a staircase of $\NN^2$. The automorphism of $k[x,y]$ exchanging $x$
and $y$ preserves the staircases. In particular it makes sense to
consider staircases in $R_i$, though the automorphism $R_i\simeq 
k[x,y]$ is not canonical. 
\\
A staircase $E\subset R_i$ defines 
a monomial zero-dimensional subcheme $Z(E)\subset U_i$ whose ideal is
generated by the monomials $m\in R_i\setminus E$.  A
multistaircase is a $r$-tuple 
$\underline E=(E_1,\dots,E_r)$ of staircases with $E_i\inc R_i$. It defines a subscheme
$Z(\underline E)=\coprod Z(E_i)$. 
\\
In our context, a $T'$-Hilbert function is a function $H:\hat{T}'\fd \NN$ such that
$\#H=\sum_{\chi \in \hat{T}'}H(\chi)$ is finite. A $T'$-Hilbert
multifunction $\underline H$ is a collection of $T'$-Hilbert functions $H_C$ parametrized
by the connected components $C$ of $X^{T'}$. Its cardinal is by
definition 
\begin{displaymath}
  \#\underline H=\sum_C \#H_C.
\end{displaymath}
Equivalently, a $T'$-Hilbert multifunction is a 
$r$-tuple  $\underline H=(H_1,\dots,H_r)$ of Hilbert functions such
that $H_i=H_{i+1}$ if $\{p_i,p_{i+1}\}\in LFix(T')$. 
\\
If $Z\subset X$ is a zero-dimensional subscheme fixed under $T'$, then $H^0(Z,\fxo_Z)$ is a
representation of $T'$ which can be decomposed as $\oplus V_{\chi}$
where $V_{\chi}\inc H^0(Z,\fxo_Z)$ is the subspace on which $T'$ acts through $\chi$. 
The $T'$-Hilbert function associated with $Z$ is by definition
$H_{T',Z}(\chi)=\dim V_\chi$. We also define a Hilbert
multifunction $\underline
H_{T',Z}$ as follows. If $p_i\in PFix(T')$, let 
$Z_i\inc Z$ the component of $Z$ located on $p_i$ and
$H_i=H_{T',Z_i}$. If $\{p_i,p_{i+1}\}\in LFix(T')$, let $Z_i=Z_{i+1}\inc Z$ the component
of $Z$ located on $V_{i,i+1}$ and $H_i=H_{i+1}=H_{T',Z_i}$. The Hilbert
multifunction associated to $Z$ is 
\begin{displaymath}
\underline
H_{T',Z}=(H_1,\dots,H_r).
\end{displaymath}
By construction, we have the equality 
\begin{displaymath}
\#\underline H_{T',Z}
=length(Z).
\end{displaymath}
A partition of $n\in \NN$ is decreasing sequence $n_1,n_2,\dots$ of
integers ( $ n_i \geq
n_{i+1}\geq 0$) with  $n_i=0$ for $i>>0$, and $\sum n_i=n$. The number of
parts is the number of integers $i$ with $n_i\neq 0$.  
\\ 
We will denote by 
\begin{listecompacte}
  \item  $Part(n)$ the
set of partitions of $n$, $Part=\coprod Part(n)$,
\item $\esc$ the set of staircases of $\NN^2$,
\item  $\mesc$ the set of multistaircases,
% \item  $\mesc(T',\underline H)$ the set of
% multistaircases $\underline E$ with $\underline H_{T',Z(\underline
%   E)}=\underline H$,
\item $\hilb(T')$ the set of $T'$-Hilbert functions,
\item $\mhilb(T')$ the set of $T'$-Hilbert multifunctions.
%\item  If $Z=\{p_{i_1},\dots,p_{i_k}\}\subset \{p_1,\dots,p_r\}$,
%  $\mesc(T',Z,H)$ is the set of parametrized by points of $W$, whose $T'$ Hilbert
%function is $H$. 
\end{listecompacte}

\subsection*{Hilbert schemes and Grassmannians}
\label{sec:hilb-scheme-grassmannians}

We denote by $\HH$ the Hilbert scheme parametrizing the 0-dimensional
subschemes of $X$. It is a disjoint union $\HH=\coprod \HH^d$, where $
\HH^d$ parametrizes the subschemes of length $d$. We denote by
$\HH_i\subset \HH$ the open subscheme parametrizing the
subschemes whose support is in $U_i$, and $\HH_{i,i+1}=\HH_i\cup
\HH_{i+1}$. 
\\ 
The action of the torus $T$ on $X$ induces an action of $T$ on
$\HH$. We denote by $\HH^T\inc \HH$ the fixed locus under this action. 
If $T'\subset T$ is a one dimensional subtorus, and $\underline H$ is
a $T'$-Hilbert multifunction, $\HH^{T',\underline H}\inc \HH$ parametrizes
by definition the subschemes $Z$, $T'$ fixed, with $T'$-Hilbert
multifunction $\underline H_{T',Z}=\underline H$. Define similarly
$\HH^{T',H}$ for a $T'$-Hilbert function $H$.  
\\
We will freely mix the above notations by intersecting the subschemes
when we gather the indexes. For instance,
$\HH_i^T=\HH_i\cap \HH^T$, $\HH^{T',H,T}=\HH^{T',H} \cap \HH^T$,
$\HH_{i,i+1}^{T'}=\HH_{i,i+1}\cap \HH^{T'}$ etc
\dots To avoid ambiguity, the formula is
\begin{displaymath}
\HH_{s}^{s_1,\dots,s_k}=\HH_{s}\cap
\HH^{s_1}\cap \dots \cap \HH^{s_k},
\end{displaymath}
where 
\begin{displaymath}
s\in \{i,\{i,i+1\}\},\ s_i\in \{d,T',(T',\underline
H),(T',H),T\}.
\end{displaymath}
\\
If $T'\hookrightarrow T$ is a one
dimensional subtorus and if  
$(i,\chi,h)\in \{1,\dots,r\}\x \hat{T}' \x \NN$, we denote by
\begin{displaymath}
\GG_{T',i,\chi,h}
\end{displaymath}
the Grassmannian parametrising the subspaces of $R_{T',i,\chi}$ of
codimension $h$. 
%We allow $h$ to be any integer with the convention
%that $\GG_{T',i,\chi,h}$ is empty if $h>\dim R_{T',i,\chi}$. 
If $H$ is a $T'$-Hilbert function, $\GG_{T',i,H}=\prod_{\chi\in \hat{T}'}
\GG_{T',i,\chi,H(\chi)}$. It is a well defined finite product 
since $G_{T',i,\chi,H(\chi)}$ 
is a point for all but
finite values of $\chi$. 
%If $H$ is a Hilbert multifunction, $\GG_{T',H}=\prod
%\GG_{T',i,H_i}$. 

\section{Description of the fixed loci}
\label{sec:descr-fixed-loci}
Let $T'\hookrightarrow T$ be a one dimensional subtorus. The goal of this section 
is to give a description of the irreducible components of $\HH^{T'}$
(proposition \ref{prop:descriptionDesComposantesIrred} and the comment
preceding it).
% \begin{prop}
%   \begin{listecompacte}\item
%     $X^T=\cup_{i}V_i$ 
%   \item  $X^{T'}=\cup_{i\in PFix(T')}V_i\cup_{(i,j)\in
%       LFix(T')}V_{ij}$.
% \end{listecompacte}
% \end{prop}
% \begin{dem}
%   The first statement is obvious. As to the second, consider the
%   decomposition of $X$ into disjoint $T$-orbits. By the very
%   definition  of a toric variety, no point of the 2-dimensional orbit
%   $T$ is 
%   $T'$-invariant. The $0$-dimensional orbits are obviously
%   $T'$-invariant. Finally, $T'$ acts on a one dimensional orbit
%   $O_{ij}=\s k[x,x^{-1}]$ by $x\mapsto t^{\alpha }x$. Thus the whole
%   orbit is fixed if $\{i,j\}\in LFix(T')$ and $O_{ij}^{T'}=\emptyset$
%   otherwise. The formula for $X^{T'}$ follows. 
% \end{dem}

\begin{thm}   \label{thr:composantesDuHilbertGradue}
\begin{displaymath}
\HH_i^{T'}=\bigcup_{H\in \hilb(T')}\HH_i^{T',H} 
\end{displaymath}
is the decomposition of
   $\HH_i^{T'}$ into smooth disjoint irreductible components. 
\end{thm}
\begin{dem}
  This is proved in \cite{evain04:irreductibiliteDesHilbertGradues}.
\end{dem}

\begin{rem}
  For some $H$,  $\HH_i^{T',H}$ may be empty so the result is that
  the irreducible components are in one-to-one correspondance with the
  set of possible Hilbert functions $H$. Throwing away the empty sets in
  the above decomposition is possible: There is an 
  algorithmic procedure to detect the emptiness of $\HH_i^{T',H}$
(\cite{evain04:irreductibiliteDesHilbertGradues}, remark 23). 
\end{rem}
Now the goal is to prove that $\HH_{i,i+1}^{T', H}$ is empty or a
product $P$ of projective spaces when $\{p_i,p_{i+1}\}\in LFix(T')$. 
An embedding $P\fd \HH^{T'}$ is
constructed in the next proposition. Then it will be shown that a non
empty $\HH_{i,i+1}^{T',H}\inc \HH^{T'}$ is the image of such an embedding. 

Let $\pi=(\pi_1,\pi_2,\dots )\in \Part(d) $ be a partition. Let
$(n_1,\dots,n_s,0=n_{s+1})$ be the finite subsequence  with no repetition 
obtained from $\pi$ with the removal of duplicates. 
Denote by
$d_l$ the number of indexes $j$ with $\pi_j=n_l$.  In other words, 
\begin{displaymath}
  \pi=(\underbrace{n_1,\dots,n_1}_{d_1\
    times},\underbrace{n_2,\dots,n_2}_{d_2\ times},\dots,
  \underbrace{n_s,\dots,n_s}_{d_s\ times},0,\dots).
\end{displaymath}
Let $\{p_i,p_{i+1}\}\in LFix(T')$ and 
$p\in V_{i,i+1}$. Since $V_{i,i+1}\subset U_i\cup U_{i+1}$, one may suppose by
  symmetry that $p\in U_i\simeq \s k[x,y]$. Exchanging the roles of
  $x$ and $y$, we may suppose that $V_{i,i+1}$ is defined by $y=0$
  in $U_i$. We denote by $Z_{p,k}$
  the subscheme with equation $(x-x(p),y^k)$. Intrinsecally, it is
  characterized as the only length $k$ curvilinear subscheme $Z\subset X$ 
  supported by $p$, $T'$-fixed, such that $Z\cap V_{i,i+1}=p$ as a
  schematic intersection.
The rational function
  \begin{eqnarray*}
    \phi_{\pi}:Sym^{d_1}V_{i,i+1}\x \dots \x Sym^{d_s}V_{i,i+1} &\dashrightarrow& \HH^{d,T'}\\
    (p_{11},\dots,p_{1d_1}),\dots,(p_{s1}\dots p_{sd_s}) &\mapsto&
    \coprod_{i\leq s,j\leq d_i} Z_{p_{ij},n_i}
  \end{eqnarray*}
is well defined on the locus where all the points $p_{a,b}\in V_{i,i+1}$ are
distinct. In fact, it is regular everywhere. 
\begin{prop}
 The function $\phi_{\pi}$ extends to
 a regular embedding  $Sym^{d_1}V_{i,i+1}\x \dots\x
 Sym^{d_s}V_{i,i+1}\fd \HH^{d,T'}$. 
\end{prop}
\begin{dem}
  The extension property is local thus it suffices to check it 
  on an open covering. The covering $V_{i,i+1}=(V_{i,i+1}\cap
  U_i)\cup (V_{i,i+1}\cap U_{i+1})=W_i\cup W_{i+1}$ 
  of $V_{i,i+1}$ induces a covering of the
  symmetric products $Sym^d V_{i,i+1}$. All the open sets in this
  covering play the same role. Thus by symmetry, it
  suffices to define an embedding 
\begin{displaymath}
\psi_\pi:Sym^{d_1}W_i\x \dots \x Sym^{d_s}
  W_i \fd \HH^{T'}
\end{displaymath}
which generically co\"\i ncides with $\phi_\pi$. 
\\
Let $Z(p_{11},\dots,p_{sd_s})\subset U_i$ be the subscheme defined by the ideal  
\begin{displaymath}
I_Z=(y^{n_1},y^{n_2}\prod_{\beta \leq
  d_1}(x-x(p_{1\beta})),\dots, y^{n_{s}}
\prod_{\alpha<s}^{\beta\leq d_{\alpha}}
  x-x(p_{\alpha\beta}),\prod_{\alpha\leq s}^{\beta\leq d_{\alpha}}
  x-x(p_{\alpha\beta})\ ).
\end{displaymath}
Let   
\begin{eqnarray*}
\psi_\pi:Sym^{d_1}W_i\x \dots \x Sym^{d_s}
  W_i &\fd& \HH^{T'}\\
_(p_{11},\dots,p_{1d_1}),\dots,(p_{s1}\dots p_{sd_s})&\mapsto& Z(p_{11},\dots,p_{sd_s}).
\end{eqnarray*}
Clearly, $Z(p_{11},\dots,p_{sd_s})$ is $T'$-fixed since $T'$ does not
act on $x$. Thus, $\psi_\pi$ is
a  well defined morphism which extends $\phi_\pi$. Now, for $1\leq \alpha
\leq s$, the transporter
$(I_Z+y^{n_{\alpha+1}+1}:y^{n_{\alpha+1}})$ 
defines a subscheme $Z_\alpha$ of $V_{i,i+1}$ of length
$d_1+\dots+d_{\alpha}$. 
Since $Z_{\alpha-1}\subset
Z_{\alpha}$ the residual scheme $Z'_{\alpha}=Z_{\alpha}\setminus
Z_{\alpha-1}$ is well defined for $\alpha \geq 2$. Consider the
morphism    
\begin{eqnarray*}
  \rho:Im(\psi_\pi) &\fd & Sym^{d_1}W_i\x \dots \x
  Sym^{d_s}W_i\\
Z&\mapsto& (Z_1,Z'_2,Z'_3,\dots,Z'_s).
\end{eqnarray*}
The composition $\rho\circ \psi_{\pi}$ is the identity. Thus,
$\psi_{\pi}$ is an embedding, as  expected. 
\end{dem}
Remark that the $T'$-Hilbert function $H_{T',Z}$ is constant when $Z$
moves in a connected
component of $\HH^{T'}$. In particular, it is constant on
$Im(\phi_{\pi})$ and $\phi_\pi$ factorizes: 
\begin{displaymath}
\phi_\pi:\prod_{\alpha \leq s} Sym^{d_\alpha}V_{i,i+1} \fd
\HH_{i,i+1}^{T',H_\pi}
\end{displaymath}
for a uniquely defined $T'$-Hilbert function $H_{\pi,T',i,i+1}$ that
we note $H_{\pi}$ for simplicity. \\
\begin{prop} \label{prop:Hii+1=produitDeProjectifs}
Let $\{p_i,p_{i+1}\}\in LFix(T')$, $H$ a $T'$-Hilbert
  function.  If  $H=H_\pi$ for some 
  $\pi\in \Part$, then $\HH_{i,i+1}^{T',H}$ is a
  product of projective spaces, thus irreducible. If $H\neq H_{\pi}$
  then $\HH_{i,i+1}^{T',H}=\emptyset$.
\end{prop}
\begin{dem}
If $H=H_\pi$ for some $\pi\in \Part$,
it suffices to prove that 
\begin{displaymath}
  \phi_\pi:\prod_{\alpha \leq s} Sym^{d_\alpha}V_{i,i+1} \fd
  \HH_{i,i+1}^{T',H_\pi}
\end{displaymath}
is an isomorphism. We already know that $\phi_{\pi}$ is an embedding thus we
   need surjectivity.  Let $Z\in \HH_{i,i+1}^{T',H_\pi}$.
   We may suppose without loss of generality that $Z\subset U_i=\s
   k[x_i,y_i]$ and that $V_{i,i+1}$ is defined by $y_i=0$ in $U_i$. 
   Since the ideal $I$ of $Z$ is
   $T'$-invariant, it is generated by elements $y_i^kP(x_i)$,
   where $P$ is a polynomial. The power $l$ being fixed, the
   polynomials $P$ such that $y_i^lP(x_i)\in I$ form an ideal in
   $k[x_i]$ generated by a polynomial $P_l$. The condition for $I$ to
   be an ideal implies the divisibility relation $P_m|P_l$ for
   $l<m$. Since $Z$ is $0$-dimensional, $P_i=1$ for $i>>0$. Let $t$ be
   the smallest integer such that $P_t=1$: 
$1=P_t|P_{t-1}|\dots|P_0$. In particular the sequence   
\begin{displaymath}
D=(D_1,D_2,\dots)=(deg(P_0),deg(P_1),\dots)
\end{displaymath}
is a partition. 
   Let $D^\nu\in \Part$ be
   the partition conjugate to $D$, ie. $D^\nu(k)=\#\{j\ s.t. \ D_j\geq
   k\}$. By construction,  $D^\nu=\pi$. 
Let $(d_1,d_2,\dots, d_s) $ be the list obtained from the list 
$(D_{t}-D_{t+1},D_{t-1}-D_{t},\dots,D_1-D_2)$ by suppression of the
zeros. Then $d_\alpha=deg(P_{j-1})-deg(P_{j})$ for some $j$
 and we let  $p_{\alpha,1},\dots,p_{\alpha,d_{\alpha}}$ be the zeros of the
   polynomial $\frac{P_{j-1}}{P_{j}}$. By definition of
   $\phi_\pi$, we have the equality
   $Z=\phi_{\pi}(p_{11},\dots,p_{sd_s})$, which shows
   the expected surjectivity.\\
If $\HH_{i,i+1}^{T',H}$ is non empty, it contains a subscheme $Z\inc X$ fixed
   under the action of $T$. Such a $Z=Z(E_i)\cup Z(E_{i+1})$ is characterized by a pair
   $(E_i,E_{i+1})$ of staircases in $R_i$ and $R_{i+1}$. Suppose as
   before that $V_{i,i+1}$ is defined by $y_i=0$ around $p_i$ and by
   $y_{i+1}=0$ around $p_{i+1}$. Using these coordinates, $E_i$ (resp
   $E_{i+1}$) is associated with a partition $\pi^i$
   (resp. $\pi^{i+1}$) defined by $x_i^ay_i^b\in E_i \Rightarrow
   b<\pi^i_{a+1}$ (resp.  $x_{i+1}^ay_{i+1}^b\in E_{i+1} \Rightarrow
   b<\pi^{i+1}_{a+1}$). Let $\pi=(\pi^i\;^\nu+\pi^{i+1}\;^\nu)^\nu$. 
   Then   
\begin{displaymath}
H=H_{T',Z}=H_{T',Z(E_i)}+H_{T',Z(E_{i+1})}=H_{\pi^i}+H_{\pi^{i+1}}=H_\pi.
\end{displaymath}
\end{dem}
Knowing that $\HH_i^{T',H}$ is empty or irreducible 
(theorem \ref{thr:composantesDuHilbertGradue}) and that
$\HH_{i,i+1}^{T',H}$ is empty or a product of projective spaces
(proposition \ref{prop:Hii+1=produitDeProjectifs}), we
obtain easily the irreducible components of $\HH^{d,T'}$:
According to the last but one item  of the next proposition,
$\HH^{T',\underline H}$ is empty or irreducible. Thus, the last item 
is the decomposition of $\HH^{d,T'}$ into irreducible components (in
fact into empty or irreducible
components and we know which terms in the union are empty).

\begin{prop}\label{prop:descriptionDesComposantesIrred}
  \begin{listecompacte}
\item $\HH^{T}=\coprod_{\underline E \in \mesc } Z(\underline E)$.
% \item     $\HH^{d,T}=\coprod_{\underline E\in \mesc, \#\underline
%     E=d}Z(\underline E)$. 
\item     $\HH^{T'}=\prod_{p_i\in PFix(T')}
\HH_i^{T'}\x \prod_{\{p_i,p_{i+1}\}\in LFix(T')}\ \HH_{i,i+1}^{T'}$. 
% \item
%   $\HH_{i,i+1}^{T'}=\coprod_{H\in \hilb}\HH^{T',H}_{i,i+1}$.
%\item     $\HH^{T'}=\coprod_{\underline H\in \mhilb(T')}\HH^{T',\underline H}$. 
\item $\HH^{T',\underline H}=\prod_{p_i\in PFix(T')}\HH_i^{T',H_i}\x
\prod_{\{p_i,p_{i+1}\}\in LFix(T')}\HH_{i,i+1}^{T',H_i}$.
\item $\HH^{d,T'}=\coprod_{\underline H\in
    \mhilb(T'),\#\underline H=d}\HH^{T',\underline H}$.
\end{listecompacte}
\end{prop}
\begin{dem}
The first point is well known. 
As to the second point, 
since the support of a subscheme $Z\inc X$
parametrised by $p\in \HH_i^{T'}$ (resp. by $p\in \HH_{i,i+1}^{T'}$)
is $p_i$ (resp. is on $V_{i,i+1}$) and since the various
$p_i,V_{i,i+1}$ do not intersect, the union morphism is a well defined
embedding 
\begin{displaymath}
\prod_{i\in PFix(T')}
\HH_i^{T'}\x \prod_{\{p_i,p_{i+1}\}\in LFix(T')}\ \HH_{i,i+1}^{T'} \fd 
\HH^{T'}.
\end{displaymath}
Since the support of a subscheme  $Z\in \HH^{T'}$ is included in
$X^{T'}=PFix(T')\cup LFix(T')$, the surjectivity is obvious. 
The third point follows from the second.
The last point is easy. 
%  Obviously, a $T$-invariant punctual subscheme $Z$ admits 
%   a $T$-invariant support. Moreover, if $p_i$ is in the support of
%   $Z$, there are natural toric coordinates around $p_i$ given by the 2
%   invariant lines through $p_i$. It is well known  that $Z$ is locally
%   defined by a monomial ideal with respect to these natural toric
%   coordinates. Thus the first point follows. The second point is
%   the first where we have taken care of the length of the subschemes.
%  The third point is again a paraphrase to say that the support of a
%   $T'$-invariant subscheme is $T'$-invariant if we note that we know
%   $X^{T'}$ according to proposition ... As to the fourth item, by
%   definition of $\HH_{ij}^{(d_1,k_1),\dots,(d_s,k_s),T'}$, there is an
%   embedding
%   $\coprod_{s,(d_1,k_1),\dots,(d_s,k_s)}\HH_{ij}^{(d_1,k_1),\dots,(d_s,k_s),T'}
%   \fd \HH_{ij}^{T'}$. Let us show that a subscheme 
% Finally, the last two points follow easily from the previous ones.
\end{dem}

\section{Equivariant Chow rings of products of Grassmannians}
\label{sec:equiv-cohom-grassm}
Let $V$ be a vector space with base $\fxb=\{e_0,\dots,e_n\}$. Let
$\chi_0,\dots,\chi_n \in \hat{T}$ be distinct characters of
$T$. These characters define an action of $T$ on $V$ by the formula 
$t.(v_0,\dots,v_n)=(\chi_0(t)v_0,\dots,\chi_n(t)v_n)$. The $T$-action 
on $V$ induces a $T$-action on the Grassmannian $\GG(d,V)$
parametrising the $d$-dimensional quotients of $V$. 
In this section, we
compute the $T$-equivariant Chow ring of $\GG(d,V)$ and of products of
such Grassmannians.

\subsection*{Equivariant Chow ring of $\GG(d,V)$}
\label{sec:equiv-cohom-ggd}
First, we recall the definition of equivariant Chow rings in the
special case of a $T$-action (To keep constant the conventions of 
the paper, we work with rational Chow groups though it is not
necessary in this section).\\ 
Let $U=(k^r\setminus 0)\x (k^r\setminus 0)$. The torus $T\simeq k^*\x
k^*$ acts on $U$
by the formula $(t_1,t_2)(v,w)=(t_1v,t_2w)$. If ${\cal X}$ is a $T$-variety, 
the quotient $(U\x {\cal X})/T=U\x^T {\cal X}$ admits a projection to $U/T=\PP^{r-1}\x \PP^{r-1}$. 
The Chow group $A_{l+2r-2}(U\x^T {\cal X})$ does not depend on the choice of
$r$ provided that $r$ is big enough (explicitly $r>\dim {\cal X}-l$)
and this Chow group is by
definition the equivariant Chow group $A_l^{T}({\cal X})$. 
In case ${\cal X}$ is smooth, we let $A^l_T({\cal X})=A_{\dim {\cal X} -l}^T({\cal X})$ and this
makes $A^*_T({\cal X})=\oplus_{l\geq 0}A^l_T({\cal X})$ a ring. 
\\
A $T$-equivariant vector bundle $F\fd {\cal X}$ defines equivariant Chern
classes: $F\x^T U$ is a vector bundle on ${\cal X}\x^T U$ and by definition
$c_i^T(F)=c_i(F\x^T U)$.  
\\
The equivariant Chow ring $A^*_T(\s k)$ of a point is a polynomial ring
by the above description. A more intrinsec description is 
as follows. 
A character $\chi\in \hat{T}$ defines canonically an equivariant 
line bundle $V_{\chi}$ over $\s k$. The map 
\begin{eqnarray*}
  \hat{T} &\fd& A^*_T(\s k)\\
  \chi &\mapsto& c_1^T(V_\chi)
\end{eqnarray*}
extends to an
isomorphism  
\begin{displaymath}
S=Sym(\hat{T}\ox \QQ)\fd A^*_T(\s k)
\end{displaymath}
where $Sym(\hat{T}\ox \QQ)$ is the symmetric $\QQ$-algebra over $\hat{T}$. \\
The morphism ${\cal X}\fd \s k$ induces by pullback a $S$-algebra structure over
$A_T^*({\cal X})$. 

Let us now turn to the case ${\cal X}=\GG(d,V)$. We denote by $\fxo(\chi)$
the line bundle $V_{\chi}\x^T U\fd U/T$. One checks easily that 
$\GG(d,V)\x^T U\fd U/T$ is the Grassmann bundle 
$\GG(d,\fxo(\chi_0)\oplus \dots \oplus \fxo(\chi_n))$. The
universal rank $d$ quotient bundle 
$Q_T\fd \GG(d,\fxo(\chi_0)\oplus \dots \oplus \fxo(\chi_n))$ 
over the Grassmann bundle 
has total space $Q_T=Q\x ^T U$ where $Q\fd \GG(d,V)$ is the
universal quotient bundle over the Grasmanniann. 
In particular $c_i^T(Q)=c_i(Q_T)$. 

If $\lambda=(\lambda_1,\dots,\lambda_{\dim V-d},0,\dots)\in \Part$, let us denote by 
\begin{displaymath}
D_{\lambda}=det(c_{\lambda_i+s-i}^T(Q))_{1\leq i,s\leq \dim V-d}
\end{displaymath}
the associated Schur
polyn{o}mial in the equivariant Chern classes of $Q$.

\begin{prop}
  The elements $D_{\lambda}$ generate the
  $S$-module $A^*_T(\GG(d,V))$.
\end{prop}
\begin{dem}
  Let $\delta\in \NN$, $A^{\leq \delta}_T \subset A^*_T(\GG(d,V))$ be the submodule defined by
  the elements of degree at most $\delta$. It suffices to prove that 
  every class in $A^{\leq \delta}_T$ is a linear combination of
  $D_{\lambda}$'s with coefficients in $S$. By definition, $A^{\leq \delta}_T =
  A^{\leq \delta}(U\x^T \GG(d,V))$ with
  $U=(k^n\setminus\{0\})\x(k^n\setminus\{0\})$ and $n>>0$. As explained, 
  the quotient $ U\x^T \GG(d,V)$ is a Grassmann bundle over $U/T$ and 
  the result follows from 
\cite{fulton84:_Intersection_theory}, Proposition 14.6.5 and 
Example 14.6.4, which describe the
  Chow ring of Grassmann bundles. 
\end{dem}

\begin{rem}\label{rem:nbGenerateursFinis}
  The number of generators in the proposition is finite since  
$D_{\lambda}=0$ for $\lambda_1>d$.
\end{rem}

To realize $A^*_T(\GG(d,V))$ as an explicit $S$-subalgebra of $S^{\GG(d,V)^T}$,
we recall from 
\cite{brion97:_equivariant_chow_groups}, 
the following result:
\begin{prop}
  If ${\cal X}$ is a projective non singular variety, the inclusion map
  $i:{\cal X}^T\fd {\cal X}$ induces an injective $S$-algebra homomorphism
  $i_T^*:A_T^*({\cal X})\fd A_T^*({\cal X}^T)$. 
\end{prop}
In fact, Brion proved the injectivity of $i_T^*$ when ${\cal X}$ is a smooth 
filtrable variety and projective varieties are filtrable. 
\\
In the present situation,
${\cal X}=\GG(d,V)$. A point $p_{\Sigma} \in \GG(d,V)^T$ is characterized by a subset
$\Sigma=\{e_{i_1},\dots,e_{i_d}\}\subset \fxb$ of cardinal $d$: if $W\subset V$
is the vector space generated by $\{e_j, e_j\notin \Sigma\}$, 
then 
\begin{displaymath}
p_{\Sigma}=V/W.
\end{displaymath}
Let $\sigma_i$ be the $i$-{th} symmetric polynomial in $d$ variables. Let 
\begin{displaymath}
c_{i,\Sigma}=\sigma_i(\chi_{i_1},\dots, \chi_{i_d})\in S
\end{displaymath}
and 
\begin{displaymath}
c_{i}\in S^{\GG(d,V)^T}=(c_{i,\Sigma})_{\Sigma\subset B,\#\Sigma=d}.
\end{displaymath}
\begin{prop} Let $i^*_T:A^*_T(\GG(d,V))\fd A^*_T(\GG(d,V)^T)=S^{\GG(d,V)^T}$ be the
  restriction morphism induced by the inclusion
  $i:\GG(d,V)^T\hookrightarrow \GG(d,V)$. Then
  $i^*_T(c_i^T(Q))=c_{i}$.
\end{prop}
\begin{dem}
The fiber of the universal quotient bundle $Q\fd \GG(d,V)$ over
$p_\Sigma$, $\Sigma=\{e_{i_1},\dots,e_{i_d}\}$, 
is a direct sum of one-dimensional representations 
with characters $\chi_{i_1},\dots,\chi_{i_d}$ thus its equivariant
total Chern class is $c^T(Q)=\prod_{j\leq d}(1+\chi_{i_j})$. 
\end{dem}

\begin{nt} If $\lambda=(\lambda_1,\dots,\lambda_{\dim V-d},0,\dots) \in \Part$, 
let 
\begin{displaymath}
\Delta_{\lambda}=det(c_{\lambda_i+s-i})_{1\leq
    i,s\leq \dim V-d}\in S^{\GG(d,V)^T}.
\end{displaymath}
As in remark \ref{rem:nbGenerateursFinis}, 
only a finite number of $\Delta_\lambda$ are non
zero. 
\end{nt}

\begin{coro}\label{prop:CohomoEquivGrassmannienneAuxPointsFixes}
  The $S$-algebra $A_T^*(\GG(d,V))\subset S^{\GG(d,V)^T}$ is generated
  as an $S$-module by the elements $\Delta_{\lambda}$. 
\end{coro}
\begin{dem}
  The restriction morphism
  $A_T^*(\GG(d,V))\fd A_T^*(\GG(d,V)^T)=S^{\GG(d,V)^T}$ is injective
  and gives the inclusion. 
  Since $c_i^T(Q)$ restricts to $c_i$, the
  generators $D_{\lambda}$ of the $S$-module  
  $A_T^*(\GG(d,V))$ restrict to $\Delta_\lambda$. The result follows. 
\end{dem}

\subsection*{Products of Grassmannians}
\label{sec:prod-grassm}

In the sequel, we will need to compute equivariant Chow rings of products of
Grassmannians, and of products in general. 
The following result explains how to deal with these
products. 

\begin{nt}
  If $P$ and $Q$ are two finite sets, $M\subset S^P$ and $N\subset S^Q$ are
  two $S$-modules, we denote by $M\ox N$ the $S$-submodule of $S^{P\x
  Q}$ image of $M\ox N$ under the natural isomorphism $ S^{P\x
  Q}\simeq S^P\ox S^Q$. 
\end{nt}

\begin{prop}\label{prop:ChowDesProduitsRestreintAuxPointsFixes}
  Let ${X}$ and $Y$ be smooth projective $T$-varieties with a finite number of fixed
  points. Let $A^*_T(X)\subset S^{X^T}$ and $A^*_T(Y)\subset S^{Y^T}$ be
  their equivariant Chow rings. Then $A^*_T(X\x Y)\subset S^{(X\x
  Y)^T}$ identifies to $A^*_T(X) \ox A^*_T(Y)$.
\end{prop}

\begin{lm}\label{lm:ChowEquivariantDuProduit}
  If $X$ and $Y$ are two smooth $T$-varieties with a finite number of fixed
  points, then $A^*_T(X\x Y)=A^*_T(X)\ox A^*_T(Y)$. 
\end{lm}
\begin{dem}
Let $F$ be a smooth variety with a cellular decomposition and $B$ a
  smooth variety. 
  According to 
\cite{edidinGraham97:CharacteristicClasses}, 
prop. 2, 
  if $\fxf\stackrel{\pi}{\fd}B$ is a locally
  trivial fibration with fiber $F$, there is a non canonical 
  isomorphism of $A^*(B)$-modules
  \begin{displaymath}
    \phi:A^*(\fxf)\fd A^*(B) \ox A^*(F).
  \end{displaymath}
  Explicitly, let us denote by
  $f_i\in A^*(F)$ the classes of the closures of the
  cells of $F$. They form a base of $A^*(F)$.  If $F_i\in A^*(\fxf)$
  is such that $F_i\cdot F=f_i$ then 
  \begin{displaymath}
    \phi^{-1}(b\ox f_i)=\pi^*b\cdot F_i.
  \end{displaymath}
In our case, $X$ and $Y$ admit cellular decompositions whose cells are
the Bialynicki-Birula strata associated with the action of a general one
parameter subgroup $T'\hookrightarrow T$. Let us denote by $V_i \subset X$ and
$W_i\subset Y$ the closures of these cells. Let 
\begin{displaymath}
X_i=V_i\x^T U \subset
X\x^T U
\end{displaymath}
and 
\begin{displaymath}
Y_i=W_i\x^T U \subset
Y\x^T U.
\end{displaymath}
By the above result about fibrations, we have:
\begin{displaymath}
  A_T^*(X)\simeq A^*(X)\ox S, \ \ A_T^*(Y)\simeq A^*(Y)\ox S,
\end{displaymath}
and the isomorphisms identify $[X_i]$ with $[V_i]\ox 1$, and  
$[Y_i]$ with $[W_i]\ox 1$.
The left and down arrows of the diagram 
\begin{eqnarray*}
  \begin{array}{ccc}
(X\x Y\x U)/T& \fd& X\x^T U\\
\fb& & \fb\\
Y\x^T U& \fd& U/T
  \end{array}
\end{eqnarray*}
yields an identification 
\begin{displaymath}
  \psi:A_T^*(X\x Y)\fd A_T^*(Y)\ox A^*(X)\fd A^*(Y)\ox S\ox A^*(X).
\end{displaymath}
Consider the natural $S$-module morphism (see
\cite{edidinGraham:constructionDesChowsEquivariants})  
\begin{displaymath}
  K:A_T^*(X)\ox_S A_T^*(Y) \fd A_T^*(X\x Y).
\end{displaymath}
The composition 
\begin{displaymath}
  \psi\circ K: A_T^*(X)\ox_S A_T^*(Y) \fd A^*(X)\ox S\ox A^*(Y)
\end{displaymath}
sends the base $[X_i]\ox [Y_j]$ to the base $[V_i]\ox 1\ox [W_j]$, thus $K$ is
an isomorphism.
\end{dem}
Now the proposition follows from the lemma and the commutativity of the following diagram.
\begin{eqnarray*}
  \begin{array}{ccc}
A_T^*(X)\ox A_T^*(Y)& \fd & A_T^*(X\x Y)\\
\downarrow i_X^*\ox i_Y^*& & \fb i_{X\x Y}^*\\
A_T^*(X^T) \ox A_T^*(Y^T)& \fd& A_T^*(X^T\x Y^T)
  \end{array}
\end{eqnarray*}

% \begin{dem}

% It suffices to prove the commutativity of the following diagram
% % By Brion (?), the equivariant Chow ring co{\"\i} ncides with the
% % equivariant cohomology. 
% % By definition,  $A^*_T(X\x Y)$ is the cohomology of the fiber product
% % $X_G \x_{\PP^{\infty}} Y_G$. Since, by Brion again (?) 
% % the cohomology of $X_G$ is a free
% % module over the cohomology of $X_G$, the Eilenberg-Moore sequence
% % degenerates and computes the cohomology of the product. Thus the
% % result.  
% \end{dem}

\section{Chow rings of graded Hilbert schemes}
\label{sec:equiv-cohom-equiv-hilbert-schemes}
Let $R=k[x,y]$. In this section, the toric variety $X$ is not projective since 
we consider the case 
$X=\s R$.
The torus $T\simeq k^*\x k^*$ acts on $X$ 
by 
\begin{displaymath}
(t_1,t_2).(x^\alpha y^\beta)=(t_1x)^\alpha (t_2y)^\beta.
\end{displaymath}
Let
$T'\hookrightarrow T$ be a one dimensional subtorus such that
$X^{T'}={(0,0)}$.  Let $H\in \hilb(T')$ be a Hilbert
function. The aim of this section is the computation of the image of
the restriction 
morphism $A_T^*(\HH^{T',H})\fd A_T^*(\HH^{T',H,T})$ (Corollary 
\ref{coro:descriptionRestrictionUsingGrassmannians}).

A point $p\in \HH^{T',H}$ parametrizes a $T'$-stable ideal 
\begin{displaymath}
I=\bigoplus_{\chi \in \hat{T'}}I_\chi\inc R,
\end{displaymath}
where $T'$ acts with character
$\chi$ on $I_{\chi}$. 
There is a $T$-equivariant embedding
\begin{eqnarray*}
\HH^{T',H}&\stackrel{l}{\hookrightarrow}&\GG_{T',H}=\prod_{\chi\in
  \hat{T}'}\GG(H(\chi),R_{\chi}) \\
I&\mapsto& (I_{\chi}). 
\end{eqnarray*}

\begin{prop}
  \begin{listecompacte}
  \item $l^*:A^*(\GG_{T',H})\fd  A^*(\HH^{T',H})$ is surjective. 
  \item $l_T^*:A^*_T(\GG_{T',H})\fd  A^*_T(\HH^{T',H})$ is surjective.
  \end{listecompacte}
\end{prop}
\begin{dem}
The surjectivity of $l^*$ has been shown by 
King and Walter
\cite{king_walter95:generateurs_anneaux_chow_espace_modules}
when $T'=\{(t,t)\}$. Their argument is valid for any $T'$ with minor
modifications. We recall briefly their method which uses ideas from 
\cite{ellingsrud-stromme93:towardsTheChowRingOfPP2}. 
Let ${\cal S}$ be an associative $k$-algebra and $M$ be a fine moduli space whose
closed points parametrize a class $\fxc$ of ${\cal S}$-modules. Denote by 
$\fxa$ the universal ${\cal S}\ox \fxo_M$-module associated with the moduli
space. 
King and Walter exhibit generators of $A^*(M)$ when 
$\fxa$ admits a nice resolution and some cohomological
conditions are satisfied.\\
In the case ${\cal S}=R$,$T'=\{(t,t)\}$, $M=\HH^{T',H}$, $\fxi$ the
universal ideal over $M$, $\fxa=(R\ox \fxo_{\HH^{T',H}} )/\fxi=\oplus \fxa_n$,
the resolution is
\begin{displaymath}
  0\fd \bigoplus_n R(-n-2)\ox_k \fxa_n \fd  \bigoplus _n R(-n-1)^2\ox_k \fxa_n
  \fd \bigoplus _n R(-n)\ox _k \fxa_n \fd \fxa \fd 0. 
\end{displaymath}
Consider now a general $T'$. For $\chi \in \hat{T}'$, we denote by
$R_{\chi}\inc R$ the subvector space on which $T'$ acts through $\chi$
and by $R(\chi)$ the $\hat T'$-graded $R$-module defined by
$R(\chi)_{\chi'}=R_{\chi+\chi'}$. Let as above $\fxa=(R\ox
\fxo_{\HH^{T',H}})/\fxi=\oplus_{\chi \in \hat{T}'}\fxa_\chi$. 
The torus $T'$ acts on $x$ and $y$
with characters $\chi_x,\chi_y$. Multiplications by $x$ and $y$ define
morphisms $\xi:\fxa_\chi\fd \fxa_{\chi+\chi_x}$ and
$\eta:\fxa_\chi\fd \fxa_{\chi+\chi_y}$. 
The resolution of $\fxa$ is:
\begin{eqnarray*}
  0\fd \bigoplus_{\chi\in\hat T'}R(-\chi-\chi_x-\chi_y)\ox_k \fxa_\chi \stackrel{\alpha}{\fd}\\
  \bigoplus _{\chi\in\hat T'} (R(-\chi-\chi_x)\oplus R(-\chi-\chi_y))\ox_k \fxa_\chi
  \stackrel{\beta}{\fd} \bigoplus _{\chi\in\hat T'} R(-\chi)\ox _k \fxa_\chi \fd \fxa \fd 0. 
\end{eqnarray*} 
where the morphisms are
\begin{displaymath}
  \alpha=\left( 
    \begin{array}{c}
-y\ox 1 +1 \ox \eta\\
x\ox 1 - 1 \ox \xi
    \end{array}
\right), \beta=(x\ox 1-1\ox \xi\ \ y\ox 1-1\ox \eta).
\end{displaymath}
With this resolution, we can follow the rest of the 
argument of
\cite{king_walter95:generateurs_anneaux_chow_espace_modules}
to conclude that $A^*(\HH^{T',H})$ is generated by the Chern classes 
$c_i(\fxa_\chi)$, hence $l^*$ is surjective.

As to the second point, remark that the morphism $l^*$ is obtained
from $l_T^*$ with the application of the functor $.\ox S/S^+$, where
$S^+\inc S$ denotes the ideal generated by the homogeneous
elements of positive degree. Since
$l^*$ is surjective, it follows from the graded Nakayama's lemma that
$l_T^*$ is surjective. 
\end{dem}
The commutative diagram
\begin{displaymath}
  \begin{array}{ccc}
\HH^{T',H}& \stackrel{l}{\hookrightarrow}& \GG_{T',H}\\
j\uparrow& & \uparrow m\\
\HH^{T',H,T}& \stackrel{n}{\hookrightarrow}& \GG_{T',H}^T
  \end{array}
\end{displaymath}
induces a map on the level of equivariant
Chow rings. Using the surjectivity of $l_T^*$, we get:
\begin{coro}\label{coro:descriptionRestrictionUsingGrassmannians}
  $Im\ j_T^*=Im\ n_T^*m_T^*$.
\end{coro}

\section{The equivariant Chow ring of $\HH^d$}
\label{sec:conclusion-proof}

Let $T'\inc T$ be a one-dimensional subtorus. 
In this section, we define finite $S$-modules 
\begin{displaymath}
{M_{T',i,H_i|}}\subset S^{\HH^{T',H_i,T}_i}
\end{displaymath}
and
\begin{displaymath}
M_{T',i,i+1,H_i} \subset S^{\HH_{i,i+1}^{T',H_i,T}}
\end{displaymath}
with explicit generators
and we prove the formula:
\begin{thm} \label{thr:description du Chow avec produit tensoriel}
  \begin{displaymath}
  A_T^*(\HH^d)=\bigcap_{T'\subset T}\bigoplus^{\underline H\in \mhilb(T')}_{
    \#\underline H=d}
  (\ \ \bigotimes_{p_i\in PFix(T')}M_{T',i,H_i|}\bigotimes_{\{p_i,p_{i+1}\}\in
  LFix(T')}M_{T',i,i+1,H_i}\ \ )
\end{displaymath}
\end{thm}
\begin{dem}
A large part of the proof consists in collecting the results from the preceding
sections using the appropriate notations. \\
Let $T'$ be a one-dimensional subtorus of $T$.
Let $p_i\in PFix(T')$. Denote by $\fxp_d(R_{T',i,\chi})$, $\chi\in \hat{T}'$ 
the set of subsets of monomials of $R_{T',i,\chi}$
%  and by $\fxp_d(R_{T',i,\chi})$ the set
% of subsets 
of cardinal $d$. A set of monomials $Z\in
\fxp_d(R_{T',i,\chi})$ defines a point $p_Z\in \GG_{T',i,\chi,d}^T$ as explained in the
preceding sections: the subspace $V_Z\inc R_{T',i,\chi}$ associated to
$p_Z$ is generated by the monomials $m\in R_{T',i,\chi}\setminus Z$. If $m\in
Z$, it is an eigenvector for the action of $T$ and we denote by
$\chi_m$ the associated character. Denote by
\begin{displaymath}
  c_{T',i,\chi,d,j,Z}=\sigma_j(\chi_m,m\in Z)\in S
\end{displaymath}
the $j$-{th} symmetric polynomial in $d$ variables evaluated on the $\chi_m$ and by 
\begin{displaymath}
  c_{T',i,\chi,d,j}=(c_{T',i,\chi,d,j,Z})_{Z\in \fxp_d(R_{T',i,\chi})}\in
  S^{\fxp_d(R_{T',i,\chi})}=S^{\GG_{T',i,\chi,d}^T}.
\end{displaymath}
For $\lambda=(\lambda_1,\dots,\lambda_{\dim R_{T',i,\chi}-d})\in
\Part$, $\lambda_1\leq d$, let 
\begin{displaymath}
  \Delta_{T',i,\chi,d,\lambda}=det( c_{T',i,\chi,d,\lambda_{r}+s-r})_{1\leq
  s,r\leq \dim R_{T',i,\chi}-d}\in S^{\GG_{T',i,\chi,d}^T}
\end{displaymath}
be the associated Schur polynomials. These Schur polynomials generate a $S$-module 
\begin{displaymath}
M_{T',i,\chi,d}\subset S^{\GG_{T',i,\chi,d}^T}
\end{displaymath}
By corollary \ref{prop:CohomoEquivGrassmannienneAuxPointsFixes}, we have
\begin{prop}
  $A_T^*(\GG_{T',i,\chi,d})\simeq M_{T',i,\chi,d}$.
\end{prop}
If $H$ is a $T'$-Hilbert function, denote 
% \begin{displaymath}
% \fxp_{T',H}(R_i)=\prod_{\chi\in
%   \hat{T}'}^{H(\chi)\neq 0}\fxp_{H(\chi)}(R_{T',i,\chi}).
% \end{displaymath}
% Let 
\begin{displaymath}
  M_{T',i,H}=\bigotimes_{H(\chi)\neq 0} M_{T',i,\chi,H(\chi)}\subset
    S^{\GG_{T',i,H}^T} 
\end{displaymath}
According to the description of equivariant Chow rings of products 
(proposition \ref{prop:ChowDesProduitsRestreintAuxPointsFixes}) and
since $\GG_{T',i,H}=\prod_{\chi\in \hat{T}'}\GG_{T',i,\chi,H(\chi)}$, we have:
  \begin{prop}
     $A_T^*(\GG_{T',i,H})\simeq M_{T',i,H}$. 
  \end{prop}
% If $\underline H=(H_1,\dots,H_n)$ is a $T'$-Hilbert multifunction and $Z\subset\{1,\dots,n\}$ parametrizes a set
% of fixed points, we denote by 
% \begin{displaymath}
% \mesc(\underline H,Z)\subset \esc^Z
% \end{displaymath}
% the subset whose 
% elements are collections of staircases $(E_i)_{i\in Z}$ with 
% $E_i\subset R_{i}$ and 
% $E_i\in \esc(T', H_i)$.
%For $Z=\{i\}$, 
The equivariant embedding 
\begin{displaymath}
\HH_i^{T',H}\hookrightarrow \GG_{T',i,H}
\end{displaymath}
yields by
restriction a morphism 
\begin{displaymath}
S^{\GG_{T',i,H}^T} \fd
S^{\HH_i^{T',H,T}}.
\end{displaymath}
If $M\subset S^{\GG_{T',i,H}^T} $, we denote by
$M_{|}$ the image of $M$ by this
restriction. 
\\
The section on the Chow ring of graded Hilbert schemes 
and corollary \ref{coro:descriptionRestrictionUsingGrassmannians}
can be reformulated in this context as:
\begin{prop}
  $A_T^*(\HH_i^{T',{H}})\simeq M_{T',i,{H}|}\subset
    S^{\HH_i^{T',H,T}}$. In particular, if $\chi_1,\dots, \chi_s \in
    \hat{T}'$ are the characters such that $H(\chi_i)\neq 0$, the generators of
 $A_T^*(\HH_i^{T',{H}})$ are the elements 
\begin{displaymath}
g_{T',i,H,\lambda_1,\dots,\lambda_s}=(\bigotimes_{\chi_j}\Delta_{T',i,\chi_j,H(\chi_j),\lambda_j})_{|}.
\end{displaymath}
\end{prop}
Now we come to the description of $A^*_T(\HH_{i,i+1}^{T',
  H})$ when $\{p_i,p_{i+1}\}\in LFix(T')$. 
Remember that we have associated a $T'$-Hilbert function $H_\pi$ to a
  partition $\pi$ such that $\HH_{i,i+1}^{T',H}\neq \emptyset$ iff
  $H=H_{\pi}$ for some $\pi$. Thus we are interested in the case
  $H=H_\pi$ and we start with the case $\pi=\pi(d,k)=(k,k,\dots,k,0,\dots)$
  where $k$ appears $d$ times. In this case, a point $p\in
  \HH_{i,i+1}^{T',H_{\pi(d,k)},T}$ parametrizes a subscheme $Z=Z_i
  \cup Z_{i+1}$ where $Z_i \in \HH_i^T$ and $Z_{i+1}\in \HH_{i+1}^T$ are
  characterized by the integers $l_{i,Z}=length(Z_i\cap V_{i,i+1})$ and
  $l_{i+1,Z}=length(Z_{i+1}\cap V_{i,i+1})=d-l_{i,Z}$ (in local
  coordinates around $p_i$ (resp. around $p_{i+1}$)
  $I_{Z_i}=(y^k,x^{l_{i,Z}})$ (resp. $I_{Z_{i+1}}=(y^k,x^{l_{i+1,Z}})$)).
\\ 
There is an action of $T$ on $V_{i,i+1}$ and we let $\chi_i$
(resp. $\chi_{i+1}=-\chi_i$) the character of $T$ which acts on the tangent space
of $p_i\in V_{i,i+1}$ (resp. of $p_{i+1}$). 
\\
For $Z\in  \HH_{i,i+1}^{T',H_{\pi(d,k)},T}$, we  define 
\begin{displaymath}
  c_{i,i+1,\pi(d,k),Z}=\frac{l_{i,Z}(l_{i,Z}+1)}{2}\chi_i+
\frac{l_{i+1,Z}(l_{i+1,Z}+1)}{2}\chi_{i+1} \in
  S.
\end{displaymath}
Then we put
\begin{displaymath}
  c_{i,i+1,\pi(d,k)}=(c_{i,i+1,\pi(d,k),Z})\in S^{ \HH_{i,i+1}^{T',H_{\pi(d,k)},T}}
\end{displaymath}
and we define 
\begin{displaymath}
  M_{T',i,i+1,H_{\pi(d,k)}}\subset S^{ \HH_{i,i+1}^{T',H_{\pi(d,k)},T}}
\end{displaymath}
to be  the $S$-module
generated by the powers  $c_{i,i+1,\pi(d,k)}^j$, $0\leq j\leq d$.
\begin{prop}
$A_T^*(\HH_{i,i+1}^{T',H_{\pi(d,k)}})\simeq
M_{T',i,i+1,H_{\pi(d,k)}}\subset S^{ \HH_{i,i+1}^{T',H_{\pi(d,k)},T}} $.
\end{prop}
\begin{dem}
  We know by proposition \ref{prop:Hii+1=produitDeProjectifs}
  that $\HH_{i,i+1}^{T',H_{\pi(d,k)}}\simeq Sym^d
  V_{i,i+1}$. Denote by $V$ the vector space with $\PP(V)=V_{i,i+1}$
  and by $P_i,P_{i+1}$ a base of $V$ with $k.P_i=p_i$,
  $k.P_{i+1}=p_{i+1}$.  
The action of $T$ on $V_{i,i+1}$ lifts to an action of $T$ on $V$ 
with characters $0$ on $P_i$ and $\chi_i$ on $P_{i+1}$. 
The action
  on $\HH_{i,i+1}^{T',H_{\pi(d,k)}}\simeq Sym^d
  V_{i,i+1}=\PP(Sym^d(V))$ is induced by  the
  characters $0,\chi_i,\dots,d\chi_i$ on $P_i^d,P_i^{d-1}P_{i+1},\dots,P_{i+1}^d$. 
Through the above identifications,  a point $Z\in \HH_{i,i+1}^{T',H_{\pi(d,k)},T}$ corresponds to the line
  $kP_i^{l_{i,Z}}P_{i+1}^{d-l_{i,Z}}\inc Sym^d(V)$. 
 In particular, the universal quotient bundle $Q=V/\fxo(-1)$
  restricts on $Z$ with equivariant Chern class  
\begin{displaymath}
c_1^TQ_{Z}=\sum_{0 \leq j\leq d,  j\neq
  d-l_{i,Z}}j.\chi_i=(\frac{d(d+1)}{2} -(d-l_{i,Z}))\chi_i.
\end{displaymath}
If we call $c_1$ the tuple $(c_1^TQ_{Z})_{Z\in {\HH_{i,i+1}^{T',H_{\pi(d,k)},T}}}$, there is a constant $a=d+1$ such that all the
coordinates of  
\begin{displaymath}
  ac_1-c_{i,i+1,\pi(d,k)}\in S^{\HH_{i,i+1}^{T',H_{\pi(d,k)},T}}
\end{displaymath}
are equal to a constant $b=-\frac{(d^2+d)}{2} \chi _i \in \ZZ\chi_i$, independent of $Z$.
The equivariant Chow ring $A_T^*(\HH_{i,i+1}^{T',H_{\pi(d,k)}} )\inc
S^{\HH_{i,i+1}^{T',H_{\pi(d,k)},T}}$ is the $S$-module generated by
the powers 
  $c_1^j$, $0\leq j\leq d$, which is also the module
  generated by the powers $c_{i,i+1,\pi(d,k)}^j$.
\end{dem}
%  \simeq Hilb^d(L_{pq})
%and that the isomorphism
%   between $\HH_{pq}^{T',H(T',pq,d,k)}$ and $Hilb^d(L_{pq})$ sends $Z$
%   to $Z\cap L_{pq}$. Thus, the computation of
%   $A_T^*(\HH_{pq}^{T',H(T',pq,d,k)})$
% and its restriction to fixed points amounts to a calculation in
%   $Hilb^d(L_{pq})\simeq \PP^d$. In a projective space, the equivariant
%   cohomology is generated by the powers of the equivariant Chern class
%   $c_1^T(\fxo(1))$, or equivalently by the powers of $c_1(\fxt)$, $\fxt$
%   being the tangent bundles, Since $c_1^T(\fxo(1))$ is a multiple a
%   $c_1^T(\fxt)$. According [mon papier a advances], 
%   the tangent bundle at
%   a point $Z\cap L_{pq}\in Hilb^d(L_{pq})^T$ is a direct sum of one
%   dimensional representations with characters
%   $\chi_p,2\chi_p,\dots,r\chi_p,\chi_q,\dots,s\chi_q$ with
%   $r=deg(Z_p\cap L_{pq})$ and $s=deg(Z_q\cap L_{pq})$. The result
%   follows.  
% \end{dem}
Let now $\pi$ be any partition and call $d_j$ then number of parts
of $\pi$ whose value is $j$. Then $H_{\pi}=\sum_{j>0}
H_{\pi(d_j,j)}$. Consider the decomposition (proposition \ref{prop:Hii+1=produitDeProjectifs})
\begin{displaymath}
  \HH_{i,i+1}^{T',H_\pi}\simeq \prod _{j>0}^{d_j>0}
 \HH_{i,i+1}^{T',H_{\pi(d_j,j)}}.
\end{displaymath}
If we adopt the convention that
\begin{eqnarray*}
  M_{T',i,i+1,H_{\pi}}&=& \bigotimes_{j>0}^{d_j>0}M_{T',i,i+1,H_{\pi(d_j,j)}}\subset 
S^{\HH_{i,i+1}^{T',H_{\pi},T} }.
\end{eqnarray*}
and if we denote by $n_1>n_2>\dots>n_s$ the integers such that
$d_{n_i}>0$, 
the formula for the equivariant Chow ring of a product yields:
\begin{prop} If $\{p_i,p_{i+1}\}\in LFix(T')$, 
   $A_T^*(\HH_{i,i+1}^{T',H_\pi})\simeq  M_{T',i,i+1,H_{\pi}} \inc
   S^{\HH_{i,i+1}^{T',H_\pi,T}}$ where $ M_{T',i,i+1,H_\pi}$ is the
   submodule generated by the elements 
\begin{displaymath}
g_{T',i,i+1,H_\pi,l_1,\dots,l_s}=\bigotimes_{j=1}^{j=s}
c_{i,i+1,\pi(d_{n_j},n_j)}^{l_j},\ \ 0\leq l_j \leq d_{n_j}.
\end{displaymath}
\end{prop}
Let $\underline H=(H_1,\dots,H_r)$ be a $T'$-Hilbert multifunction. The decomposition
\begin{displaymath}
  \HH^{T',\underline H}\simeq \prod_{p_i\in PFix(T')}\HH_i^{T',H_i}\prod_{\{p_i,p_{i+1}\} \in
  LFix(T')} \HH_{i,i+1}^{T',H_i}
\end{displaymath}
yields the following formula for $A_T^*(\HH^{T',H})$. 
\begin{prop}
  $A_T^*(\HH^{T',\underline H})\subset S^{\HH^{T',\underline H,T}}$ identifies to the
$S$-module 
\begin{displaymath}
\bigotimes_{p_i\in PFix(T')}M_{T',i,H_i|}\bigotimes_{\{p_i,p_{i+1}\}\in
  LFix(T')}M_{T',i,i+1,H_i} .
\end{displaymath}
\end{prop}
Since 
\begin{displaymath}
\HH^{d,T'}=\coprod_{\underline H\in
    \mhilb(T'),\#\underline H=d}\HH^{T',\underline H},
\end{displaymath}
we obtain:
\begin{prop}
   $A_T^*(\HH^{d,T'})\subset S^{\HH^{d,T}}$ is isomorphic to
\begin{displaymath}
 \bigoplus_{\underline H\in
    \mhilb(T'),\#\underline H=d\ }(\bigotimes_{\ p_i\in PFix(T')}M_{T',i,H_i|}\bigotimes_{\{p_i,p_{i+1}\}\in
  LFix(T')}M_{T',i,i+1,H_i}\ )
\end{displaymath}
\end{prop}
Now by \cite{brion97:_equivariant_chow_groups}, theorem 3.3,
\begin{displaymath}
A_T^*(\HH^d)\simeq \bigcap_{T'\subset
  T}A_T^*(\HH^{d,T'}),
\end{displaymath}
which proves the theorem.
\end{dem}

\section{Description by congruences}
\label{sec:descr-congr}
In the last section, 
$A_T^*(\HH^{d})$ has been described 
by a formula involving tensor products
and intersections. The goal of this section is to give a simpler
presentation. Explicitly, we will 
describe $A_T^*(\HH^{d})\inc S^{\HH^{d,T}}$
as a set of tuples of elements of $S$ satisfying congruence
relations. 
 
The possibility to reformulate the description of the last section
with congruence relations was suggested to me by Michel Brion. 

Let $\pi:{\cal X}\fd \s k$ be a smooth projective $T$-variety with a
finite number of 
fixpoints, and
\begin{displaymath}
\begin{array}{ccccc}
  f:A_T^*({\cal X})&\ox_S& A_T^*({\cal X}) &\fd & S=A_T^*(\s k)\\
x&\ox& y &\mapsto & \pi_*(x\cdot y).
\end{array}
\end{displaymath}
Let $Q=Frac(S)$. According to the localisation theorem 
(\cite{brion97:_equivariant_chow_groups}, cor. 3.2.1)
the morphism
$i_T^*:A_T^*({\cal X}) \hookrightarrow  S^{{\cal X}^T}$ becomes an isomorphism
\begin{displaymath}
i^*_{T,Q}:A_T^*({\cal X})_Q=A_T^*({\cal X})\ox_S Q \fd Q^{{\cal X}^T} 
\end{displaymath}
after tensorisation with $Q$. % In particular $f_Q:A_T^*({\cal X})_Q \ox_S&
% A_T^*({\cal X})_Q &\fd & Q$
% is well defined on $S^{{\cal X}^T}\ox S^{{\cal X}^T}$.

\begin{prop}
  Let $\beta_i=i_{T,Q}^*(\overline \beta_i)$ be a set of generators of the $S$-module
  $i_T^*A_T^*({\cal X})\subset S^{{\cal X}^T}$ and $\alpha=i_{T,Q}^*(\overline \alpha) \in S^{{\cal X}^T}$. Then
  $\alpha \in  i_{T}^*A_T^*({\cal X}) \Leftrightarrow \forall i, f_Q(\overline \alpha\ox
  \overline \beta_i)\in S$.
\end{prop}
\begin{dem}
  According to lemma \ref{lm:appartenanceAUnReseauParProduitScalaire}
  applied with $M=i_T^*A_T^*({\cal X})$, 
  it suffices to find bases $B_i, C_j \in
  A_T^*({\cal X})$ with $f(B_i\ox C_j)=\delta_{ij}$. 
  Let $\lambda:T'=k^*\subset T$ be a one
parameter subgroup with ${\cal X}^{T'}={\cal X}^T$.
The Bialynicki-Birula cell associated to a point $p_i\in
  {\cal X}^{T'}$ is the set $\{x\in {\cal X},\ \lim_{t\fd \infty}\lambda(t).x=p_i\}$. We
  denote by $b_i^+$ its closure and we let $B_i^+=b_i^+\x^T U\subset {\cal X}\x^T U$. It
follows from the proof of lemma \ref{lm:ChowEquivariantDuProduit} 
that the elements $[B_i^+]\in A_T^*({\cal X})$ form a
$S$-base. Consider similarly the cells $B_i^-$ defined by the one
  parameter subgroup $\lambda\circ i$, where $i:k^*\fd k^*,\ x\mapsto
  x^{-1}$. 
% Order (partially) the points $p_1,\dots,p_k \in {\cal X}^{T'}$ by the
%   relation $p_i<p_j$ if there is an irreducible curve $C$, $T'$-stable,
%   such that for a general point $p\in C$, $\lim_{t\fd \infty}\lambda(t).p=p_j$
%   and $\lim_{t\fd 0}\lambda(t).p=p_i$. By construction, 
It is a property of the Bialynicki-Birula cells that one can order the
points $p_i$ such that:
  \begin{eqnarray*}
    b_i^+ \cap b_j^- &\neq& \emptyset \Rightarrow p_j\leq p_i,\\
    b_i^+ \cap b_i^- &=&p_i \ \ \  (\mathrm{transversal\ intersection}).
  \end{eqnarray*}
It follows that 
 \begin{eqnarray*}
    f([B_i^+]\ox [ B_j^-]) &\neq& 0 \Rightarrow p_j\leq p_i,\\
    f([B_i^+] \ox [B_i^-] &=& 1.
  \end{eqnarray*}
% We choose a total order on the $p_i$'s which refines the partial order
% introduced above. 
Up to relabelling, one may suppose $p_1<p_2<\dots
<p_n$. The matrix $m_{ij}=f([B_i^+]\ox [ B_j^-])$ is a lower triangular
unipotent matrix. In particular, there exists a triangular matrix
$\lambda_{ij}$ such that $[C_j]=\sum \lambda_{ij}[B_i^-]$ verifies
$f([B_i^+]\ox [C_j])=\delta_{ij}$.
\end{dem}

\begin{lm}\label{lm:appartenanceAUnReseauParProduitScalaire}
  Let $Q=Frac(S)$, $M\subset S^n$ be a free $S$-module, $M_Q=M\ox_S Q$, $f:M\ox_S
  M\fd S$ be $S$-linear, $f_Q:M_Q\ox M_Q \fd Q$ be the  $Q$-linear map
  extending $f$, and $\beta_1,\dots,\beta_p$ be generators of $M$. Suppose
  that 
  \begin{listecompacte}
    \item $M\inc S^n$ yields an isomorphism $M_Q\simeq Q^n$ after
    tensorisation with $Q$,
    \item there exist basis $(B_1,\dots,B_n)$, $(C_1,\dots,C_n)$ of
      $M$ such that $f(B_i\ox C_j)=\delta_{ij}$.
  \end{listecompacte}
Let $\alpha \in S^n$. Then $\alpha\in M \Leftrightarrow \forall i, f_Q(\alpha\ox
  \beta_i)\in S$. \findem
\end{lm}
%  \begin{dem}
%    The implication $\Rightarrow$ is obvious. As to the implication
%    $\Leftarrow$, one may suppose by linearity that
%    $\beta_i=C_i$. Write $\alpha \in S^n\ox Q=M\ox Q$ as
%    $\sum \alpha_i B_i$, $\alpha_i\in Q$. Since
%    $\alpha_i=f(\alpha\ox C_i)\in S$, $\alpha \in M$. 
%  \end{dem}

As a corollary, we get a description of $i_T^*A_T^*({\cal X})\inc S^{{\cal X}^T}$ in terms of
congruences involving equivariant Chern classes of the
restrictions $T_{{\cal X},p}$ of the tangent bundle $T_{\cal X}$ to fixed points. 
\begin{coro}
  Let $\beta_i=(\beta_{ip})_{p\in {\cal X}^T}$ be a set of generators of the $S$-module
  $i_T^*A_T^*({\cal X})\subset S^{{\cal X}^T}$ and $\alpha=(\alpha_p)\in S^{{\cal X}^T}$. Then the
  following conditions are equivalent.
  \begin{listecompacte}
  \item 
    $\alpha \in i_T^*A_T^*({\cal X}) $
  \item $\forall i,\ \sum_{p \in {\cal X}^T} (\alpha_p\beta_{ip} \prod_{q\neq
  p}c_{\dim {\cal X}}^T(T_{{\cal X},q}))\equiv 0\ (\prod_{p \in {\cal X}^T}c_{\dim {\cal X}}^T(T_{{\cal X},p}))$
\end{listecompacte}
\end{coro}
\begin{dem}
Let us write $\beta_i=i_{T,Q}^*(\overline \beta_i)$, $\alpha=i_{T,Q}^*(\overline \alpha)$.
By the integration formula of Edidin and Graham 
\cite{edidin_Graham98:formuleDeBott},  $f_Q(\overline \alpha\ox
  \overline \beta_i)=  \sum_{p \in {\cal X}^T}
  \frac{\alpha_p\beta_{ip}}{c_{\dim {\cal X}}^T(T_{{\cal X},p})}$.
Thus, the corollary is nothing but the criteria of the last
proposition. 
\end{dem}

We can collect in a set 
\begin{math}
G(T',\underline{H})\inc S^{\HH^{T',\underline H,T}}
\end{math}
the generators of
$i_T^*A_T^*(\HH^{T',\underline{H}}) \inc  S^{\HH^{T',\underline{H},T}}$
constructed in section \ref{sec:conclusion-proof}.
Explicitly, $G(T',\underline H)$ contains
the elements 
\begin{displaymath}
  g_{T',\underline H,\lambda_{ij},l_{ij}}=\bigotimes _{p_i\in PFix(T')}g_{T',i,H_i,\lambda_{i1},\dots,\lambda_{i,s_i}}\bigotimes_{\{p_i,p_{i+1}\}\in LFix(T')} g_{T',i,i+1,H_i,l_{i1},\dots,l_{i,t_i}}.
\end{displaymath}
These generators and the last corollary make
it possible to obtain a description of $i_T^*A_T^*(\HH^{T',\underline{H}})$ via
congruences. To get a description of 
\begin{displaymath}
A_T^*(\HH^{d})\simeq \bigcap_{T'\inc
  T\ }\bigoplus_{\ \HH^{T',\underline{H}}\neq \emptyset}  i_T^*A_T^*(\HH^{T',\underline{H}})
\end{displaymath}
we merely have to gather the congruence relations constructed for
the various $\HH^{T',\underline{H}}$. We finally obtain:
\begin{thm} \label{thr:description du Chow avec congruences}
  The ring $A_T^*(\HH^{d})\inc S^{\HH^{d,T}}$ is the set of tuples
  $\alpha=(\alpha_{p})$ such that, 
  $\forall T'\inc T$ one-dimensional subtorus , $\forall  \underline{H}\in
  \mhilb(T')$ with $\HH^{T',\underline{H}}\neq \emptyset$, 
  $\forall g=(g_{p})\in G(T',\underline H)$,
the congruence relation    
\begin{displaymath}
 \sum_{p \in {\HH^{T',\underline{H},T}}} (\alpha_{p}g_{p} \prod_{q\neq
   p}^{q \in {\HH^{T',\underline{H},T}}}
c_{\dim {\HH^{T',\underline{H}}}}^T(T_{{\HH^{T',\underline{H}}},q}))\equiv 0\ (\prod_{p \in {\HH^{T',\underline{H},T}}}c_{\dim
  {\HH^{T',\underline{H}}}}^T(T_{{\HH^{T',\underline{H}}},p}))
\end{displaymath}
holds. 
\end{thm}

\begin{rem}
  The tangent space at a $T$-fixed point of
  $\HH_i^{T',H}$ or $\HH_{i,i+1}^{T',H}$  is
  known \cite{evain04:irreductibiliteDesHilbertGradues}. In
  particular, since $\HH^{T',\underline{H}}$ is a product of terms
  isomorphic to  $\HH_i^{T',H}$ or $\HH_{i,i+1}^{T',H}$, the
  equivariant Chern
  classes appearing in the theorem are explicitly computable (See
  the example in the next section).
\end{rem}

\begin{rem}
 
\end{rem}

Let $S^+=\hat TS\subset S$. 
The description of the usual Chow ring now follows from 
\cite{brion97:_equivariant_chow_groups},cor.2.3.1.
\begin{thm}
  The ring $A^*(\HH^d)$ is the quotient of $A_T^*(\HH^d)\subset
  S^{\HH^{d,T}}$ 
  by the ideal $S^+A_T^*(\HH^{d})$ generated by the elements
  $(f,\dots,f)$, $f\in S^+$.
\end{thm}

\begin{rem} For $d=1$, $\HH^d=X$ and one recovers that 
  $A_T^*(X)$ is isomorphic to the space of continuous piecewise
  polynomial functions on the fan of $X$.

  The Betti numbers of $\HH^d$ can be computed using the description
  of the last two
  theorems. In particular, one can check for small values of
  $d$ and explicit surfaces $X$ that the Betti numbers are those
  computed with G\"ottsche's formula.  
\end{rem}

\section{An example}
\label{sec:an-example}

In this section, we compute the Chow ring of the Hilbert
scheme $\HH^3=\HH^3 \plp$.

First, we fix the notations: $T=k^*\x k^*=\s k[t_1^{\pm 1},t_2^{\pm
  1}]$  and $\plp=Proj\ k[x_1,x_2,x_3]$. 
The torus $T$ acts on $\plp$ and on itself. The symmetric group $S_3$
acts on $\plp$. The action of an element $(a,b)\in T$, $\sigma \in S_3$ 
is as follows. 
\begin{eqnarray*}
  (a,b).x_1^\alpha x_2 ^\beta x_3^\gamma=x_1^\alpha (ax_2)
  ^\beta (bx_3)^\gamma\\
  (a,b).t_1^\alpha t_2^{\beta}=(at_1)^\alpha (bt_2)^{\beta}\\
   \sigma.x_1^\alpha x_2 ^\beta x_3^\gamma=x_{\sigma(1)}^\alpha x_{\sigma(2)} ^\beta x_{\sigma(3)}^\gamma.
\end{eqnarray*}
The equivariant map $T\fd \plp$, $(a,b)\fd (1,a,b)$ identifies $t_1$
with $\frac{x_2}{x_1} $, and $t_2$ with $\frac{x_3}{x_1} $.
We denote by $p_1=(1:0:0),p_2=(0:1:0),p_3=(0:0:1)$ the  three toric points
of $\plp$. The plane $\plp$ is covered by the three affine planes
$U_1=\s k[t_1,t_2]= \s R_1$, $U_2=\s k[t_1^{-1},t_1^{-1}t_2]=\s R_2$, $U_3=\s
k[t^{-1}_2, t_1t_2^{-1}]=\s R_3$. 
Since $S_3$ acts on $T=\{x_1x_2x_3\neq 0\}$, 
it acts on $\hat{T}$ by
$\sigma.\chi(t)=\chi(\sigma^{-1}t)$, and on $S=Sym(\hat{T}\ox \QQ)$. 
If $T'\inc T$ and $H\in \hilb(T')$ is a
$T'$-Hilbert function, let $\sigma.H\in \hilb(\sigma.T')$ be the Hilbert
function defined by $(\sigma.H)(\chi)=H(\sigma^{-1}.\chi)$. 
If $\underline H=(H_1,H_2,H_3)\in \mhilb(T')$ is a Hilbert multifonction, let
$\sigma.\underline H\in \mhilb(\sigma.T')$ be the Hilbert multifunction with
$(\sigma.\underline H)_i=\sigma.H_j$ where $j$ is such that
$\sigma.p_j=p_i$. 
To each subvariety $\HH^{T',\underline H}\inc \HH$, we have associated
a set of congruence relations $R_i$. Explicitly, constants $u_i\in S$ and $d_i(q)\in
S$ for $q\in \HH^{T',\underline
  H,T}$ have been defined 
such that $s \in
S^{\HH^{T',\underline H,T}}$ satisfies $R_i$ if
\begin{displaymath}
  \sum_{q\in \HH^{T',\underline H,T}} d_i(q)s(q)  \equiv 0(u_i).
\end{displaymath}
The subvariety $\sigma.\HH^{T',\underline
  H}=\HH^{\sigma.T',\sigma.\underline H}$ is
associated with the set of congruence relations $\sigma.R_i$ where by
definition $s\in \HH^{\sigma.T',\sigma\underline H,T}$ satisfies
$\sigma.R_i$ if :
\begin{displaymath}
 \sum_{q\in \HH^{\sigma.T',\sigma.\underline H,T}} d_i(\sigma^{-1}q)s(q)  \equiv 0(\sigma.u_i).
\end{displaymath}
Summing up, there is an action of $S_3$ on the set of congruence
relations. We will produce the set of relations up to this action. 

We list the possible $p\in \HH^T$. Let $E_1=\{1,t_1,t_2\}\inc R_1$,
$E_2=\{1,t_1,t_1^2\}\inc R_1$, $E_3= \{1,t_1,\}\inc R_1$, $E_4= \{1\}\inc R_1$, 
$E_5=\{1\}\inc R_2$, $E_6=\{1\}\inc R_3$.
The multistaircases 
\begin{eqnarray*}
  \underline E_A=(E_1,\emptyset,\emptyset)\\
  \underline E_B=(E_2,\emptyset,\emptyset)\\
    \underline E_C=(E_3,E_5,\emptyset)\\
  \underline E_D=(E_3,\emptyset,E_6)\\
  \underline E_E=(E_4,E_5,E_6)
\end{eqnarray*}
are associated with points $A,B,C,D,E \in \HH^T$. Up to the action of
$S^3$, these are the only points of $\HH^{3,T}$.  
  \begin{figure}[h] 
     \begin{center}
        \input{multiescaliers.pstex_t}
     \end{center} 
  \end{figure}

We recall the description of the tangent space at $p\in \HH^T$ where
$p$ is described by a multistaircase $(F_1,F_2,F_3)$ 
(\cite{evain04:irreductibiliteDesHilbertGradues}). The staircase
$F_i$ is a set of monomials in $R_i=k[x,y]$ where $x,y$ are the
toric coordinates around $p_i$. 
A cleft for $F_i$ is a monomial $m=x^ay^b \notin F_i$ with ($a=0$ or
$x^{a-1}y^b\in F_i$) and ($b=0$ or
$x^{a}y^{b-1}\in F_i$). We order the clefts of $F_i$ according to their
$x$-coordinates: $c_1=y^{b_1},c_2=x^{a_2}y^{b_2},\dots,c_p=x^{a_p}$ 
with $a_1=0<a_2<\dots<a_p$. 
An $x$-cleft couple for $F_i$ is a couple $C=(c_k,m)$, where $c_k$ is a cleft
($k\neq p$),
$m\in F_i$, and $mx^{a_{k+1}-a_k}\notin F_i$. The torus $T$ acts on the
monomials $c_k$ and $m$ with characters $\chi_k$ and $\chi_m$. We
let $\chi_C=\chi_m-\chi_k$. By symmetry, there is a notion of
$y$-cleft couple for $F_i$. The set of cleft
couples for $p$ is by definition the union of the ($x$ or $y$)-cleft
couples for $F_1$, $F_2$, $F_3$. The vector space $T_p\HH$ is in
bijection with the formal sums $\sum \lambda_i C_i$, where $C_i$ is a
cleft couple for $p$. Moreover, under this correspondance, the cleft couple
$C$ is an eigenvector for the action of $T$ 
with character $\chi_C$. 

If $p\in \HH^{T}$, and if $\underline H$ is the $T'$-Hilbert multifunction of the
subscheme associated with $p$, we let $\HH^{T',p}=\HH^{T',\underline H}$.    
The subvariety $\HH^{T',p}\inc \HH$ give non trivial congruences only if
$\HH^{T',p}$ is not a point, ie. if  
$T_p\HH^{T',p}\neq 0$. Using the above description of the tangent
space, we find for each point $p$ a finite number of possible $T'$.
The results are collected in
the following array. Under each point $p$ are listed the couples
$(a,b)$ such that $T_{ab}=\{t^a,t^b\}\inc T$ verifies
$\HH^{T_{ab},p}\neq \{p\}$. For each such $(a,b)$, the
corresponding dimension $\dim \HH^{T_{ab},p}$ is given. 
\begin{displaymath}
  \begin{array}{|c|c|c|c|c|}
\hline
A& B& C& D& E\\
\hline
\begin{array}{cc}
a,b& dim\\
1,0& 2\\
0,1& 2\\
2,1& 1\\
1,2& 1
\end{array}
&
\begin{array}{cc}
a,b& dim\\
1,0& 1\\
0,1& 3\\
1,1& 1\\
1,2& 1
\end{array}
&
\begin{array}{cc}
a,b& dim\\
1,0& 1\\
0,1& 3\\
1,1& 2
\end{array}
&
\begin{array}{cc}
a,b& dim\\
1,0& 2\\
0,1& 2\\
1,1& 2
\end{array}
& 
\begin{array}{cc}
a,b& dim\\
1,0& 2\\
0,1& 2\\
1,1& 2
\end{array}\\
\hline
 \end{array}
\end{displaymath}
For some $a,b,p$, $a',b',p'$, we have an identification
$\HH^{T_{ab},p}=\sigma.\HH^{T_{a'b'},p'}$ ($\sigma\in  S^3$). 
Explicitly, up to action, we have 
$H^{T_{10},A}=H^{T_{10},D}=H^{T_{01},A}$, 
$H^{T_{01},B}=H^{T_{01},C}$, 
$H^{T_{12},A}=H^{T_{12},B}=H^{T_{21},A}$, 
$H^{T_{11},C}=H^{T_{11},D}$, 
$H^{T_{01},D}=H^{T_{01},E}=H^{T_{10},E}=H^{T_{11},E}$.
Thus, by symmetry, we only consider $(a,b,p)$
within the following list:
\begin{displaymath}
  \{(0,1,A),(1,2,A),(1,0,B),(0,1,B),(1,1,B),(1,1,C),(1,0,C),(0,1,D)\}.
\end{displaymath}
For each of the above values of $(a,b,p)$, we construct the congruence
relations associated with the variety $\HH^{T_{ab},p}$. The results
are summed up in the following array.
\\
\epsfig{file=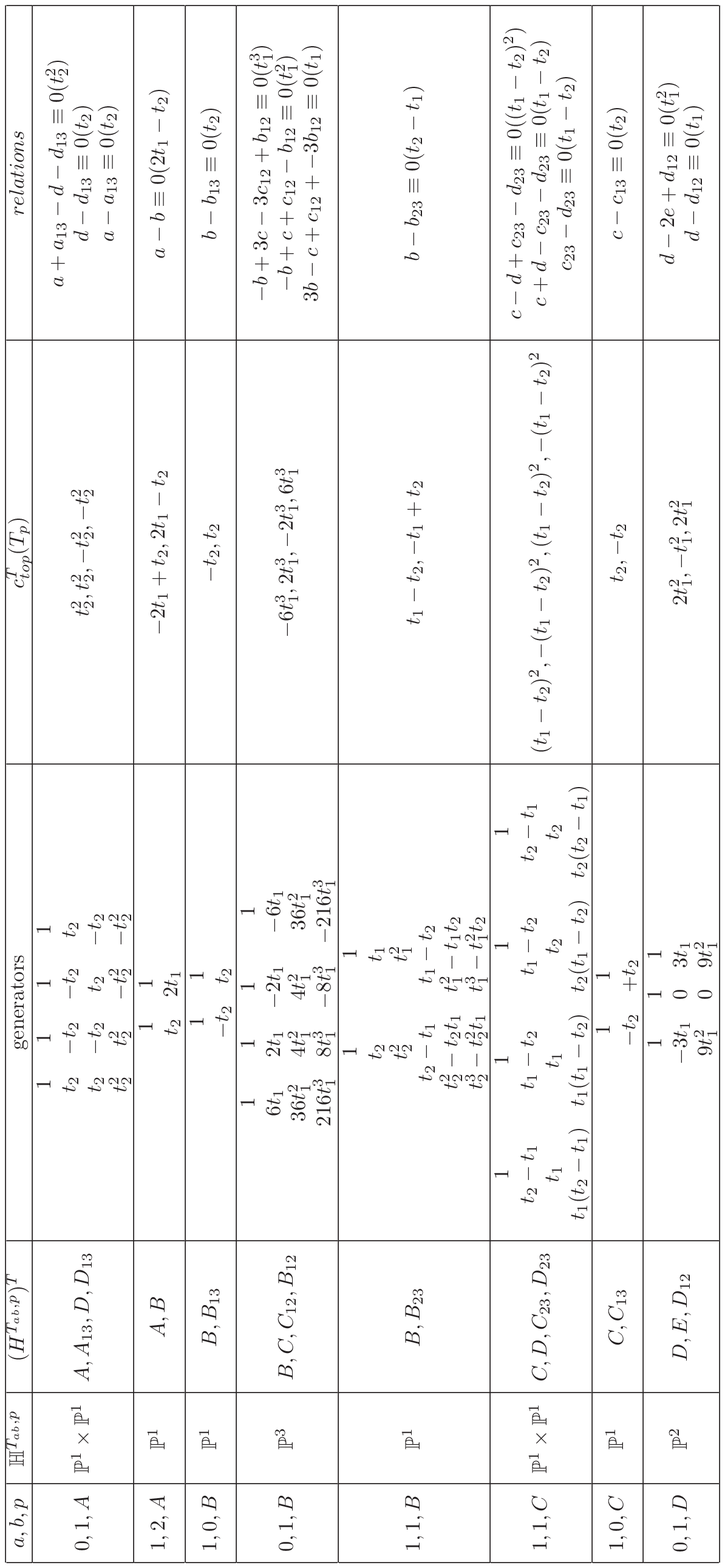,height=220mm,width=170mm,angle=0}
If $\sigma=(n_1,n_2)\in S_3$ is a permutation and $p\in \HH^T$, we
denote by $p_{n_1n_2}$ the element $\sigma.p$. We explain
how to read the array, taking the second line as an example. 
The first three column mean that
$\HH^{T_{01},A}$ is isomorphic to $\PP^1\x \PP^1$ and contains the
points $A,A_{13},D,D_{13}$.  Four generators of
$A_T^*(\HH^{T_{01},A})\inc S^{\{A,A_{13},D,D_{13}\}}$ have been constructed in
section \ref{sec:conclusion-proof}, namely
$A+A_{13}+D+D_{13}$,\dots,$t_2^2A+t_2^2A_{13}-t_2^2D-t_2^2D_{13}$. 
The coefficients of these expressions are written down in the fourth column. 
The top equivariant Chern classes  $c_{top}^T(T_{A}\HH^{T_{01},A})$, \dots, 
$c_{top}^T(T_{D_{13}}\HH^{T_{01},A})$ are respectivly $t_2^2$, \dots,
$-t_2^2$,
as indicated in the fifth column. 
We can construct congruence relations with these
data following the procedure of section \ref{sec:descr-congr}:
$A_T^*(\HH^{T_{01},A})\inc S^{\{A,A_{13},D,D_{13}\} }$ is the
set of elements $aA+a_{13}A_{13}+dD+d_{13}D_{13}$ whose
coefficients $a,\dots,d_{13}$ verify $a+a_{13}-d-d_{13}\equiv 0(t_2^2)$, 
$d-d_{13}\equiv 0(t_2)$ and
$a-a_{13}\equiv 0(t_2)$. This is the meaning of the last column. 
We gather the congruence relations constructed in the array, and we obtain:
\begin{thm}\label{thr:leCasHilbTroisP2}
  The equivariant Chow ring  $A_T^*(\HH^3\plp)\inc
  \QQ[t_1,t_2]^{\{A,A_{12},\dots,E\}}$ 
is the set of linear
  combinations $aA+a_{12}A_{12}+\dots +eE$ satisfying the relations
  \begin{listecompacte}
%1ere ligne de donnees
\item 
$a+a_{13}-d-d_{13}\equiv 0(t_2^2)$
\item $d-d_{13}\equiv 0(t_2)$
\item $a-a_{13}\equiv 0(t_2)$
%2eme ligne de donnees
\item
$a-b\equiv 0(2t_1-t_2)$
%3eme ligne de donnees
\item
$b-b_{13}\equiv 0(t_2)$
%4eme ligne de donnnees
\item
$-b+3c-3c_{12}+b_{12}\equiv 0(t_1^3)$
\item $-b+c+c_{12}-b_{12}\equiv 0(t_1^2)$
\item $3b-c+c_{12}+-3b_{12}\equiv 0(t_1)$
%5eme ligne de donnees
\item
$b-b_{23}\equiv 0(t_2-t_1)$
%6eme ligne de donnees
\item
 $ c-d+c_{23}-d_{23} \equiv 0((t_1-t_2)^2)$
\item   $c+d-c_{23}-d_{23} \equiv 0(t_1-t_2)$
\item $c_{23}-d_{23} \equiv 0(t_1-t_2)$
%7eme ligne de donnees
\item
$c-c_{13}\equiv 0(t_2)$
%8eme ligne de donnees
\item
$d-2e+d_{12}\equiv 0(t_1^2)$
\item
$d-d_{12}\equiv 0(t_1)$
\item all relations deduced from the above by the action of the
  symmetric group $S_3$. 
  \end{listecompacte}
\item 
The Chow ring $A^*(\HH^3\plp)$ is the quotient of $A_T^*(\HH^3\plp)$
by the ideal generated by the elements $fA+\dots +fE$, $f\in \QQ[t_1,t_2]^+$. 
\end{thm}

\bibliography{/users/evain/perso/fichiersConfig/modeles/texBiblio.bib}
\bibliographystyle{plain} 

\end{document}

%% file: multiescaliers.pstex_t
\begin{picture}(0,0)%
\includegraphics{multiescaliers.pstex}%
\end{picture}%
\setlength{\unitlength}{1243sp}%
\begingroup\makeatletter\ifx\SetFigFont\undefined%
\gdef\SetFigFont#1#2#3#4#5{%
  \reset@font\fontsize{#1}{#2pt}%
  \fontfamily{#3}\fontseries{#4}\fontshape{#5}%
  \selectfont}%
\fi\endgroup%
\begin{picture}(17577,3394)(661,-3059)
\put(676,-2986){\makebox(0,0)[lb]{\smash{{\SetFigFont{5}{6.0}{\rmdefault}{\mddefault}{\updefault}{\color[rgb]{0,0,0}$p_1$}%
}}}}
\put(3601,-2761){\makebox(0,0)[lb]{\smash{{\SetFigFont{5}{6.0}{\rmdefault}{\mddefault}{\updefault}{\color[rgb]{0,0,0}$p_2$}%
}}}}
\put(676,164){\makebox(0,0)[lb]{\smash{{\SetFigFont{5}{6.0}{\rmdefault}{\mddefault}{\updefault}{\color[rgb]{0,0,0}$p_3$}%
}}}}
\put(2026,-1861){\makebox(0,0)[lb]{\smash{{\SetFigFont{5}{6.0}{\rmdefault}{\mddefault}{\updefault}{\color[rgb]{0,0,0}A}%
}}}}
\put(5401,-1861){\makebox(0,0)[lb]{\smash{{\SetFigFont{5}{6.0}{\rmdefault}{\mddefault}{\updefault}{\color[rgb]{0,0,0}B}%
}}}}
\put(8776,-1861){\makebox(0,0)[lb]{\smash{{\SetFigFont{5}{6.0}{\rmdefault}{\mddefault}{\updefault}{\color[rgb]{0,0,0}C}%
}}}}
\put(12376,-1861){\makebox(0,0)[lb]{\smash{{\SetFigFont{5}{6.0}{\rmdefault}{\mddefault}{\updefault}{\color[rgb]{0,0,0}D}%
}}}}
\put(15751,-1861){\makebox(0,0)[lb]{\smash{{\SetFigFont{5}{6.0}{\rmdefault}{\mddefault}{\updefault}{\color[rgb]{0,0,0}E}%
}}}}
\end{picture}%

%% file: chow.bbl
\begin{thebibliography}{10}

\bibitem{brion97:_equivariant_chow_groups}
M~Brion.
\newblock Equivariant chow groups for torus actions.
\newblock {\em Transformations Groups}, (2):225--267, 1997.

\bibitem{costello-grojnowski03:CohomoSchemaHilbertPonctuel}
K~Costello and I~Grojnowski.
\newblock Hilbert schemes, hecke algebras and the calogero-sutherland system.
\newblock {\em math.AG/0310189}.

\bibitem{edidinGraham97:CharacteristicClasses}
Dan Edidin and William Graham.
\newblock Characteristic classes in the {C}how ring.
\newblock {\em J. Algebraic Geom.}, 6(3):431--443, 1997.

\bibitem{edidinGraham:constructionDesChowsEquivariants}
Dan Edidin and William Graham.
\newblock Equivariant intersection theory.
\newblock {\em Invent. Math.}, 131(3):595--634, 1998.

\bibitem{edidin_Graham98:formuleDeBott}
Dan Edidin and William Graham.
\newblock Localization in equivariant intersection theory and the {B}ott
  residue formula.
\newblock {\em Amer. J. Math.}, 120(3):619--636, 1998.

\bibitem{ellingsrud-stromme87:chow_group_of_hilbert_schemes}
Geir Ellingsrud and Stein~Arild Str{\o}mme.
\newblock On the homology of the {H}ilbert scheme of points in the plane.
\newblock {\em Invent. math.}, (87):343--352, 1987.

\bibitem{ellingsrud-stromme93:towardsTheChowRingOfPP2}
Geir Ellingsrud and Stein~Arild Str{\o}mme.
\newblock Towards the {C}how ring of the {H}ilbert scheme of {${\bf P}\sp 2$}.
\newblock {\em J. Reine Angew. Math.}, 441:33--44, 1993.

\bibitem{evain04:irreductibiliteDesHilbertGradues}
Laurent Evain.
\newblock Irreducible components of the equivariant punctual {H}ilbert schemes.
\newblock {\em Adv. Math.}, 185(2):328--346, 2004.

\bibitem{fantechi-gottsche93:cohomologie-3-points}
Barbara Fantechi and Lothar G{\"o}ttsche.
\newblock The cohomology ring of the {H}ilbert scheme of {$3$} points on a
  smooth projective variety.
\newblock {\em J. Reine Angew. Math.}, 439:147--158, 1993.

\bibitem{fulton84:_Intersection_theory}
W.~Fulton.
\newblock {\em Intersection Theory}.
\newblock Number~3 in Ergebnisse der Mathematik. Springer, 2 edition, 1984.

\bibitem{gottsche90:nbBettiDuSchemaHilbertDesSurfaces}
Lothar G{\"o}ttsche.
\newblock Betti numbers for the {H}ilbert function strata of the punctual
  {H}ilbert scheme in two variables.
\newblock {\em Manuscripta Math.}, 66(3):253--259, 1990.

\bibitem{grojnowski:resultatsSimilairesANakajima}
I.~Grojnowski.
\newblock Instantons and affine algebras. {I}. {T}he {H}ilbert scheme and
  vertex operators.
\newblock {\em Math. Res. Lett.}, 3(2):275--291, 1996.

\bibitem{haiman_sturmfels02:multigradedHilbertSchemes}
Mark Haiman and Bernd Sturmfels.
\newblock Multigraded {H}ilbert schemes.
\newblock {\em J. Algebraic Geom.}, 13(4):725--769, 2004.

\bibitem{king_walter95:generateurs_anneaux_chow_espace_modules}
A~King and C~Walter.
\newblock On chow rings of fine moduli spaces of modules.
\newblock {\em J. reine angew. Math.}, 461:179--187, 1995.

\bibitem{lehn99:_chern_classes_of_tautological_sheaves_on_Hilbert_schemes}
M~Lehn.
\newblock Chern classes of tautological sheaves on hilbert schemes of points on
  surfaces.
\newblock {\em Invent. math.}, (136):157--207, 1999.

\bibitem{lehn_sorger02:cup_product_on_Hilbert_schemes_for_K3}
Manfred Lehn and Christoph Sorger.
\newblock The cup product of the hilbert scheme for k3 surfaces.
\newblock {\em Invent. Math (to appear) and [math.AG/0012166]}.

\bibitem{lehn_sorger01:cup_product_on_Hilbert_schemes}
Manfred Lehn and Christoph Sorger.
\newblock Symmetric groups and the cup product on the cohomology of {H}ilbert
  schemes.
\newblock {\em Duke Math. J.}, 110(2):345--357, 2001.

\bibitem{liQinWang04mathAG:cohomoDesSchemaHilbertSurface=FibreSurP1}
Wei-Ping Li, Zhenbo Qin, and Weiqiang Wang.
\newblock The cohomology rings of {H}ilbert schemes via {J}ack polynomials.
\newblock In {\em Algebraic structures and moduli spaces}, volume~38 of {\em
  CRM Proc. Lecture Notes}, pages 249--258. Amer. Math. Soc., Providence, RI,
  2004.

\bibitem{nakajima97:_heisenberg_et_Hilbert_schemes}
H~Nakajima.
\newblock Heisenberg algebra and {H}ilbert schemes of points on projective
  surfaces.
\newblock {\em Ann. of Math.}, 2(145):379--388, 1997.

\bibitem{vasserot01:anneauCohomologieHilbert}
Eric Vasserot.
\newblock Sur l'anneau de cohomologie du sch{\'e}ma de {H}ilbert de {$\bold
  C\sp 2$}.
\newblock {\em C. R. Acad. Sci. Paris S{\'e}r. I Math.}, 332(1):7--12, 2001.

\bibitem{Vezzosi_Vistoli:KTheorieEquivarianteEtSousTores}
Gabriele Vezzosi and Angelo Vistoli.
\newblock Higher algebraic {$K$}-theory for actions of diagonalizable groups.
\newblock {\em Invent. Math.}, 153(1):1--44, 2003.

\end{thebibliography}
